\documentclass[12pt]{report}
\usepackage{amssymb}

\newtheorem{dfn}{Definition}
\newtheorem{lmm}{Lemma}
\newtheorem{rem}{Remark}

\newtheorem{thrm}{Theorem}
\newenvironment{proof}{\par\smallskip\noindent{\bf Proof:}}
                      {\hfill $\blacksquare$\\[1.5ex]}
\newcommand\trace[1]{\mbox{\rm trace}(#1)}
\newcommand\ov{\overline}
\newcommand\val[1]{\bigm| _{#1}}
\newcommand\sm{\setminus}
\newcommand\fa{\forall}

\renewcommand\d{\mbox{\rm d}}

\newcommand\res[2]{\mbox{\rm Res}(#1;#2)}
\renewcommand\det[1]{|#1|}

\newcommand\der[2]{\frac{\d{#1}}{\d{#2}}}
\newcommand\dder[2]{\frac{\partial{#1}}{\partial{#2}}}
\newcommand\Set[1]{\{#1\}}
\newcommand\ie{{\it i.e.\ }}

\newcommand\diag[1]{{\mbox{\rm diag}}(#1)}
\newcommand\rank[1]{\mbox{\rm rank}(#1)}
\newcommand\Ker[1]{\mbox{\rm Ker}(#1)}
\renewcommand\Im[1]{\mbox{\rm Im}(#1)}
\title{ The Schlesinger System and the
Riemann-Hilbert Problem} \author{{Dan Volok}} \date{}
\begin{document} \maketitle \begin{abstract}
We generalize some
classical results for the Schlesinger system of partial
differential equations and give the explicit form of its
solution, associated with rational matrix functions in general
position.
\end{abstract} \tableofcontents
\chapter{Introduction}
 The Schlesinger system first appeared in
L.~Schlesinger's work \cite{schl} as a completely integrable non-
linear Pfaffian system, governing the isomonodromic deformations
in the class of non-resonant Fuchsian systems. It is closely
related with the Riemann-Hilbert monodromy problem, which requires
to find a Fuchsian system with prescribed monodromy. This problem
was included by D.~Hilbert in his  list of problems \cite{hilb1}
as the  21st problem and solved  in  the non-resonant case  by D.~
Hilbert himself in \cite{hilb3}, I.~Plemelj in \cite{plem} and
G.~D.~Birkhoff in \cite{bir}. These solutions are based on passing
from the Riemann-Hilbert monodromy problem to the Riemann-Hilbert
boundary problem, which requires to find the factorization of a
matrix function on a contour into the product of two factors,
holomorphic and invertible, respectively, inside and outside the
contour. In \cite{miwa} T.~Miwa has employed a new method, based
on the holonomic quantum fields theory and Clifford operators, to
construct the isomonodromic deformation of a non- resonant
Fuchsian system. Assuming the appropriate restrictions on the
initial values, he has proved that all the singularities, except
the fixed ones, of a solution of the Schlesinger system are poles.

The main goal of
this Thesis is to investigate
the Schlesinger system without such restrictions.
Omitting the assumption of non-resonance, we encounter several
difficulties.
First of all,   the class of isomonodromic
deformations of a resonant Fuchsian system can be very rich. An
example of this phenomenon, concerning rational matrix functions
in general position,  was described by V.~E.~Katsnelson in
\cite{katze1}, \cite{katze}.
Although the Schlesinger
system preserves monodromy in general,
it is unclear {\it a priori}, how to distinguish the isomonodromic
deformation, corresponding to its solution.
On the other hand,
A.~A.~Bolibruch has demonstrated (see \cite{anbo}) that in general
the Riemann-Hilbert monodromy problem may have no solution at all.
 We cannot overcome these difficulties by adapting
 T.~Miwa's construction  because of
its heavy reliance on the  assumption of non-resonance.
Instead, we achieve our goal
in the following way.

Given a linear differential system with rational coefficients (not
necessarily Fuchsian), we can factorize its fundamental solution
in a neighborhood of  a singular point into the product of two
factors. One of them is holomorphic and invertible in a
neighborhood of  the singular point, and we call it the
non-singular factor. The other is a "translation" of a matrix
function, holomorphic and invertible in the Riemann surface of the
logarithm. We  call this matrix function the principal factor.
 The existence of such local factorization was established
by G.~D.~Birkhoff (see \cite{bir1}), using his method
of solution of the   Riemann-Hilbert boundary problem.
Furthermore, we consider
 isomonodromic\footnote{The monodromy being
understood in the strictly "Fuchsian" sense, without taking into
account  the Stokes phenomena.}  deformations, which
preserve also all the principal factors.
We  call such deformations
 isoprincipal. The isoprincipal deformation of a given linear differential
system is uniquely determined, and
 we derive a
completely integrable non-linear Pfaffian
system, which governs the isoprincipal deformations.
We  call it the generalized Schlesinger system.
In the case of Fuchsian systems this is just the classical
Schlesinger system, which explains the terminology.
The Cauchy problem
for the generalized Schlesinger system can be reduced
to the Riemann-Hilbert boundary problem, depending on a parameter.
We solve this problem by D.~Hilbert's
method, based on  Fredholm theory of linear integral
equations, and obtain a generalization of Miwa theorem for
the generalized Schlesinger system.

The presentation of our results is organized as follows.
\begin{itemize}
\item  Chapter~\ref{intro} is dedicated
to G.~D.~Birkhoff's and D.~Hilbert's methods for the Riemann-Hilbert
boundary problem. We present these classical results in some
detail, since they are going to play the crucial role in our
derivations.
\item  In Chapter~\ref{Schlesinger} we
define the generalized Schlesinger system and the isoprincipal deformation and
formulate two of the main results of this Thesis.
 These are Theorem~\ref{schlthrm}, which states that the isoprincipal
deformations
 are  governed by the generalized
Schlesinger system, and Theorem~\ref{tm}, which states
that all the singularities, except the fixed ones,
of a solution of the generalized Schlesinger system are poles.
\item In Chapter~\ref{Fuchs} we apply our results, obtained in
Chapter~\ref{Schlesinger}, to the case of Fuchsian systems.
In particular, we show that Schlesinger theorem
(Theorem~\ref{Schthrm} here), formulated for the non-resonant Fuchsian
systems, is a special case of our Theorem~\ref{schlthrm}.
 We also consider V.~E.~Katsnelson's example, concerning
rational matrix functions in general position, and formulate
the third main result of this Thesis -  Theorem~\ref{explicit}.
It gives the  explicit form of solution of the Schlesinger system in
this case. This result was obtained in joint work with V.~E.~Katsnelson.
\end{itemize}
We would also like to mention, in addition to
the literature cited above, the following textbooks:
\begin{itemize} \item on  general theory of linear differential
systems with rational coefficients: \cite{hille}, \cite{gol},
\cite{codlev},  \cite{sib};
\item on the Schlesinger system:
\cite{iksy};
\item on the Riemann-Hilbert boundary problem:
\cite{hilb2}, \cite{gakh}, \cite{ablo}. \end{itemize}
\chapter{The Riemann-Hilbert Boundary Problem} \label{intro}
\section{A Theorem of G. D. Birkhoff}
 Let $\gamma$ be a closed analytic Jordan
curve   in the complex plane ${\mathbb C},$ let ${\cal D}$ be the
simply connected domain in ${\mathbb C},$ bounded by $\gamma:$
\begin{equation} \gamma=\partial{\cal D}, \end{equation}
 and let
 ${\cal A}$ be an open neighborhood  of $\gamma.$
We denote by ${\cal D}^+,{\cal D}^-$ the  domains
\begin{eqnarray}
\label{d+}
{\cal D}^+&=&{\cal D}\cup{\cal A},\\
\label{d-} {\cal D}^-&=&\Set{{\mathbb C}\sm{\cal D}}\cup{\cal A},
\end{eqnarray}
so that
\begin{eqnarray}
\label{cov1}
{\cal D}^+\cup{\cal D}^-&=&{\mathbb C},\\
\label{cov2} {\cal D}^+\cap{\cal D}^-&=&{\cal A}.
\end{eqnarray}
We also denote by ${\cal O}({\cal A})$ the algebra of all
functions,
 holomorphic\footnote{By "holomorphic" we mean "analytic univalued".} in ${\cal A},$
 by ${\cal O}^*({\cal A})$ the group of all functions,
 holomorphic and nowhere vanishing in ${\cal A},$
  by $gl(m,{\cal O}({\cal A}))$ the algebra of all $m\times m$ matrix functions,
holomorphic in ${\cal A},$ and by $GL(m,{\cal O}({\cal A}))$ the
group of all $m\times m$ matrix functions,
 holomorphic and invertible in ${\cal A}.$
\begin{dfn}
Let $F(x)\in GL(m,{\cal O}({\cal A})).$ The Riemann-Hilbert
boundary problem (for $F(x)$ on $\gamma$) is to find matrix
functions
\begin{eqnarray} X^+(x)&\in&GL(m,{\cal O}({\cal D}^+)), \label{mulan+} \\
X^-(x)&\in&GL(m,{\cal O}({\cal D}^-)), \label{mulan-}
\end{eqnarray}
such that \begin{equation}\label{rhp} X^-(x)F(x)=X^+(x),\quad
x\in{\cal A}.
\end{equation}
\end{dfn}
\begin{rem}
We would like to note here  that
\begin{enumerate}
\item ${\cal A}$ can always be replaced in (\ref{mulan+}) -
(\ref{rhp}) by a smaller neighborhood of $\gamma,$ because
$X^+(x)$ and $X^-(x)$ admit the analytic continuation wherever $Q$
does; \item $X^-(x),$  in general, need not be holomorphic or
invertible at $\infty;$ \item $X^+(x)$ and $X^-(x)$ are determined
up to multiplication by an arbitrary entire invertible matrix
function from the left.
\end{enumerate}
\end{rem}
Let us start dealing with the Riemann-Hilbert boundary
problem by recalling the following classical result,
which follows immediately from Cauchy and Liouville theorems:
\begin{lmm}
\label{cauchy} Let $Z(x)\in gl(m,{\cal O}({\cal A})).$ Then there
exist unique
\begin{eqnarray}
Z^+(x)&\in&gl(m,{\cal O}({\cal D}^+)),\label{addan+} \\
Z^-(x)&\in&gl(m,{\cal O}({\cal
D}^-\cup\Set{\infty})),\label{addan-}
\end{eqnarray}
such that
\begin{equation}\label{addpr}
 Z(x)=Z^+(x)-Z^-(x),\quad x\in{\cal A}
\end{equation}
and
\begin{equation}\label{addpr0}
Z^-(\infty)=0.
\end{equation}
These matrix functions $Z^+(x), Z^-(x)$ are given by
\begin{eqnarray}
\label{dz+}
 Z^+(x)&=&\left\{\begin{array}{l@{\quad}l}
\displaystyle\oint_{\gamma}\frac{Z(y)}{y-x}\frac{\d y}{2\pi i},& x\in{\cal D}, \\
\displaystyle Z(x)+ \oint_{\gamma}\frac{Z(y)-Z(x)}{y-x}\frac{\d
y}{2\pi i},& x\in{\cal A},
\end{array}\right.\\
\label{dz-}
Z^-(x)&=&\left\{\begin{array}{l@{\quad}l}
\displaystyle\oint_{\gamma}\frac{Z(y)}{y-x}\frac{\d y}{2\pi i},&
x\in {\mathbb C}\sm\ov{\cal D},\\
 \displaystyle\oint_{\gamma} \frac{Z(y)-Z(x)}{y-x}\frac{\d y}{2\pi i},&x\in{\cal A}.
\end{array}\right.
\end{eqnarray}
\end{lmm}
Let us observe that Lemma~\ref{cauchy} provides a solution to the
Riemann-Hilbert boundary problem in the case $m=1,$ \ie when
$F(x)\in{\cal O}^*({\cal A}).$
 Indeed, let $\mu\in{\mathbb Z}$ be the index of $F(x):$
\begin{equation}\mu=\oint_{\gamma}\der{F}{x}F^{-1}(x)\frac{\d x}{2\pi i}.
\end{equation} Then, having fixed a point $x_0\in {\cal D}\sm{\cal
A},$ we can define the {\em univalued} function \begin{equation}
Z(x)=\log(F(x)(x-x_0)^{-\mu})\in{\cal O}({\cal A}). \end{equation}
In view of Lemma~\ref{cauchy}, we can let functions
$Z^+(x),Z^-(x)$  satisfy relations (\ref{addan+}) - (\ref{addpr0})
for $m=1$ and define
\begin{eqnarray}
X^+(x)&=&\exp(Z^+(x)), \\
 X^-(x)&=&(x-x_0)^{-\mu}\exp(Z^-(x)).
\end{eqnarray}
Then we can conclude that
relations (\ref{mulan+}) - (\ref{rhp}) hold true. Thus the Riemann-Hilbert
boundary problem  in this case has  solution $X^+(x),X^-(x),$
such that $X^-(x)$ has at most a pole at $\infty.$
Although this method is unsuitable for  $m>1,$ it turns out that the
same conclusion is also true in the general case.
This result is due to G.~D.~Birkhoff (see \cite{bir}):
\begin{thrm}[Birkhoff] \label{Birthrm}
Let $F(x)\in GL(m,{\cal O}({\cal A})).$ Then there exist $X^+(x),
X^-(x),$ such that relations~(\ref{mulan+}) - (\ref{rhp}) are
satisfied and $X^-(x)$ has at most a pole at $\infty.$ \end{thrm}
\begin{proof}
We shall present  a sketch of the  proof, which can be found in
\cite{bir}. For simplicity, we assume \begin{equation} {\cal
D}=\Set{x:|x|<1}. \end{equation} In the general case, the
reasonings below can be slightly modified, using  the Riemann
conformal map theorem.

The main idea is to consider  the
following system of  boundary problems:
\begin{equation}\label{birsys}\left\{ \begin{array}{r@{=}l}
Z_1^+(x)-Z_1^-(x)&x^\nu Z_2^+(x)F(x), \\
Z_2^+(x)-Z_2^-(x)&I-x^{-\nu}Z_1^-(x)F^{-1}(x),\\
Z_1^-(\infty)=Z_2^-(\infty)&0,
 \end{array} \right. \end{equation}
where $\nu$ is a positive integer.
We shall show that, for $\nu$
sufficiently large , system~(\ref{birsys}) is solvable by
successive approximations.

Indeed, let us look for solution in the form
\begin{equation}\label{serz}
Z^\varsigma_j(x)=\sum_{k=0}^{\infty}Z_{j,k}^\varsigma(x), \quad
j=1,2,\varsigma=+,-, \end{equation} where $Z_{j,k}^\varsigma(x)$
are recursively defined by \begin{equation} \left\{
\begin{array}{r@{=}l}
Z_{1,k+1}^+(x)-Z_{1,k+1}^-(x)&x^\nu Z_{2,k}^+(x)F(x), \\
Z_{2,k+1}^+(x)-Z_{2,k+1}^-(x)&-x^{-\nu}Z_{1,k}^-(x)F^{-1}(x), \\
Z_{1,k}^-(\infty)=Z_{2,k}^-(\infty)&0
\end{array} \right. \end{equation}
with the initial conditions
\begin{equation} Z_{1,0}^+(x)\equiv 0, Z_{1,0}^-(x)\equiv 0,
 Z_{2,0}^+(x)\equiv I, Z_{2,0}^-(x)\equiv 0. \end{equation}
Furthermore, let $\epsilon\in(0,1)$ be such that
\begin{equation}\Set{x:\epsilon\leq|x|\leq\frac{1}{\epsilon}}\subset{\cal A},
\end{equation} and for $Z\in gl(m,{\cal O}({\cal A}))$ let us denote
\begin{equation}\|Z\|=\max_{{1\leq\alpha,\beta\leq m\atop
\epsilon\leq|x|\leq\frac{1}{\epsilon}}} |Z(x)_{\alpha\beta}|.
\end{equation} Now for $\epsilon<|x|<\frac{1}{\epsilon}$ \begin{equation}
\begin{array}{ll} \displaystyle
Z_{1,k+1}^-(x)=&\displaystyle\oint_{|y|=1 } \frac{y^\nu
Z_{2,k}^+(y)F(y)-x^\nu Z_{2,k}^+(x)F(x)}{y-x}
\frac{\d y}{2\pi i}= \\
&\displaystyle \oint_{|y|=\epsilon }
\frac{y^\nu Z_{2,k}^+(y)F(y)}{y-x} \frac{\d y}{2\pi i}= \\
&\displaystyle\oint_{|y|=\epsilon }
 y^\nu Z_{2,k}^+(y)\frac{F(y)-F(x)}{y-x}\frac{\d y}{2\pi i}.
\end{array} \end{equation}
Analogously,
\begin{equation}
Z_{2,k+1}^+(x)= -\oint_{|y|=\frac{1}{\epsilon }}
y^{-\nu}Z_{1,k}^-(y)\frac{F^{-1}(y)-F^{-1}(x)}{y-x} \frac{\d y}{2\pi i}.
\end{equation}
Hence  we obtain
\begin{eqnarray}
\|Z_{1,k+1}^-\|&\leq&\epsilon^\nu a \|Z_{2,k}^+\|,\\
\|Z_{2,k+1}^+\|&\leq&\epsilon^\nu a \|Z_{1,k}^-\|,
\end{eqnarray}
where \begin{equation} a=\max\Set{\max_{{1\leq\alpha,\beta\leq
m\atop \epsilon\leq|x|,|y|\leq\frac{1}{\epsilon}}}
\left|\frac{F(x)_{\alpha\beta}-F(y)_{\alpha\beta}}{x-y}\right|,
\max_{{1\leq\alpha,\beta\leq m\atop
\epsilon\leq|x|,|y|\leq\frac{1}{\epsilon}}} \left
|\frac{F^{-1}(x)_{\alpha\beta}-F^{-1}(y)_{\alpha\beta}}{x-y}\right
|},\end{equation} and uniform convergence of series (\ref{serz})
for $\nu$ sufficiently large follows immediately.

Now we can define
\begin{eqnarray}
X^+_0(x)&=&Z_1^+(x), \\
X^-_0(x)&=&x^{\nu}(I+Z_2^-(x)),
\end{eqnarray}
and observe that
\begin{eqnarray}
X^+_0(x)&\in&gl(m,{\cal O}({\cal D}^+)), \\
X^-_0(x&\in&gl(m,{\cal O}({\cal D}^-)), \\
X^-_0(x)F(x)&=&X^+_0(x), \quad x\in{\cal A},
\end{eqnarray}
and $X^-_0(x)$ has at most a pole at $\infty.$ Also, since
$Z_2^-(\infty)=0,$ $\det{X^-_0(x)}$ does not identically vanish,
and, therefore, neither does $\det{X^+_0(x)}.$ The isolated zeroes
of $\det{X^+_0(x)},\det{X^-_0(x)}$ can now be successively
eliminated by the following procedure. Let $x_0\in{\cal D}$ be a
zero of $\det{X^+_0(x)},$ then there exists $T\in GL(m,{\mathbb
C}),$ such that the first row of $TX^+_0(x)$ vanishes at $x_0$ (if
$x_0\in{\cal A}$ then the same is true for $TX^-_0(x)$). Let us
consider
\begin{eqnarray}
X_1^+(x)&=&E_{x_0}(x)TX_0^+(x),\\
X_1^-(x)&=&E_{x_0}(x)TX_0^-(x),
\end{eqnarray}
 where  $E_{x_0}(x)=\diag{\frac{1}{x-x_0},1,\ldots,1}.$
Then
\begin{eqnarray}
X^+_1(x)&\in&gl(m,{\cal O}({\cal D}^+)), \\
X^-_1(x&\in&gl(m,{\cal O}({\cal D}^-)), \\
X^-_1(x)F(x)&=&X^+_1(x), \quad x\in{\cal A},
\end{eqnarray}
  and the multiplicity of
zero of $\det{X^+_1(x)}$ at $x_0$ is less by 1.
 Repeating this procedure, we
finally obtain $X^+(x),X^-(x),$ satisfying~(\ref{mulan+}) -
(\ref{rhp}), where $X^-(x)$  has at most a pole at $\infty.$
\end{proof}
One can also consider the Riemann-Hilbert boundary problem when
${\cal D}$ is a finite disjoint union of simply connected domains,
rather than a single domain. Let us denote for the moment
\begin{equation} {\cal D}=\bigcup_{j=1}^n{\cal D}_{j},
\end{equation}
 where for $j=1,\ldots,n$  ${\cal D}_j$ is
a bounded   simply connected domain in ${\mathbb C},$ such that
\begin{equation}\ov{\cal D}_{j^{\prime}}\cap\ov{\cal
D}_{j^{\prime\prime}}=\emptyset,\quad j^{\prime}\not=
j^{\prime\prime}. \end{equation} Then $\gamma=\partial{\cal
D}=\cup\gamma_j,$ where $\gamma_j=\partial{\cal D}_{j},$ and we
can assume that neighborhood ${\cal A}$ of $\gamma$ is a disjoint
union of neighborhoods ${\cal A}_j$ of $\gamma_j:$
\begin{eqnarray}
{\cal A}&=&\bigcup_{j=1}^n{\cal A}_j, \\
 {\cal A}_{j^{\prime}}&\cap&{\cal A}_{j^{\prime\prime}}=\emptyset,
 \quad j^{\prime}\not= j^{\prime\prime}.
\end{eqnarray}
Solution of the  Riemann-Hilbert
boundary problem in this case  can be reduced to the previously
considered one.  Indeed, as we have already shown, for
$j=1,\ldots,n$ one  can successively find
\begin{eqnarray}
X_{j}^+(x)&\in& GL(m,{\cal O}({\cal D}_{j}\cup{\cal A}_j)), \\
X_{j}^-(x)&\in&GL(m,{\cal O}(\Set{{\mathbb C}\sm{\cal
D}_j}\cup{\cal A}_j)),
 \end{eqnarray}
such that \begin{equation} X_{j}^-(x)\left (X_{j-1}^-(x)\cdots
X_1^-(x)F(x)\right )= X_{j}^+(x), \quad x\in{\cal A}_{j}.
\end{equation} Define
\begin{eqnarray}
X^+(x)&=&X_n^-(x)\cdots X_{j+1}^-(x)X_{j}^+(x), \quad x\in{\cal D}_{j}\cup{\cal A}_j, \\
X^-(x)&=&X_n^-(x)\cdots X_1^-(x),\quad x\in{\cal D}^-,
\end{eqnarray}
then $X^+,X^-$ give a solution to the Riemann-Hilbert boundary problem
for $F(x).$ If each $X_{j}^-(x)$ has at most
a pole at $\infty$ then so does $X^-(x),$ \ie the conclusion
of Theorem~\ref{Birthrm} holds true also in this case.

We would also like to mention another approach to the
Riemann-Hilbert boundary problem, concerning holomorphic vector
bundles. In view of (\ref{cov1}), (\ref{cov2}) matrix function
$F(x)\in GL(m,{\cal O}({\cal A}))$ can be used to define
 a  holomorphic vector bundle over ${\mathbb C}.$ The existence theorem
of H.~Grauert (see \cite{grau1}) states that any holomorphic
vector bundle over an open Riemann surface (in particular, over
${\mathbb C}$) is holomorphically isomorphic to the trivial one.
This means that there exist
\begin{eqnarray}
X^+(x)&\in&GL(m,{\cal O}({\cal D}^+)), \\
X^-(x)&\in&GL(m,{\cal O}({\cal D}^-)),
\end{eqnarray}
such that \begin{equation} X^-(x)F(x)=X^+(x),\quad x\in{\cal
A}.\end{equation} Thus the Riemann-Hilbert boundary problem is
always solvable. Furthermore, we can use $F(x)$ to construct a
holomorphic vector bundle over $\ov{\mathbb C}.$ It was proved
 by A.~Grothendieck in \cite{gro} that any holomorphic vector
bundle over the Riemann sphere is holomorphically
isomorphic to a direct sum of line bundles.
Recalling our solution of  the Riemann-Hilbert boundary problem
in the case $m=1$ (the case of a line bundle), we can conclude
that there exist
\begin{eqnarray} X^+(x)&\in&GL(m,{\cal O}({\cal D}^+)), \\
X^-(x)&\in&GL(m,{\cal O}({\cal D}^-\cup\Set{\infty})),
\end{eqnarray}
and $\mu_1,\ldots,\mu_m\in{\mathbb Z},$ such that
\begin{equation}\label{birgro1} X^-(x)F(x)=D(x)X^+(x),\quad
x\in{\cal A}, \end{equation} where
\begin{equation} D(x)=\diag{(x-x_0)^{mu_1},\ldots,(x-x_0)^{\mu_m}},
\end{equation} $x_0\in{\cal D}\sm{\cal A}.$ This is, of course, a
stronger result than
 Theorem~\ref{Birthrm}.
A similar result was obtained by
 G.~D.~Birkhoff in \cite{bir2}. It states that there
exist
\begin{eqnarray} X^+(x)&\in&GL(m,{\cal O}({\cal D}^+)), \\
X^-(x)&\in&GL(m,{\cal O}({\cal D}^-\cup\Set{\infty})),
\end{eqnarray}
and $\mu_1,\ldots,\mu_m\in{\mathbb Z},$ such that
\begin{equation}\label{birgro2} X^-(x)F(x)=X^+(x)D(x),\quad
x\in{\cal A}, \end{equation} where
\begin{equation} D(x)=\diag{(x-x_0)^{mu_1},\ldots,(x-x_0)^{\mu_m}},
\end{equation} $x_0\in{\cal D}\sm{\cal A}.$ It should be noted,
however, that in general factors
 $X^+,X^-,D$
are different in (\ref{birgro1}) and (\ref{birgro2}).
\section{Regular Solutions, Depending on a
Parameter}
Let us investigate, when the Riemann-Hilbert
boundary problem has  solution
$X^+(x), X^-(x),$ such that $X^-(x)$ is  holomorphic and invertible also at
$\infty.$ We shall call such a solution {\em regular.}
First of all, we observe that if $X^+(x), X^-(x)$ give a regular solution
to   the Riemann-Hilbert boundary problem for $F(x),$
then the determinants of $X^+(x), X^-(x)$ give a
regular solution to   the Riemann-Hilbert boundary problem for
the determinant of $F(x):$
\begin{eqnarray}
\det{X^+(x)}&\in&{\cal O}^*({\cal D}^+), \\
\det{X^-(x)}&\in&{\cal O}^*({\cal D}^-\cup\Set{\infty}), \\
\det{X^-(x)}\det{F(x)}&=&\det{X^+(x)},\quad x\in{\cal A}.
\end{eqnarray} Hence a necessary condition for the existence of  a
regular solution is \begin{equation}\label{conrhp}
\mu=\oint_{\gamma} \der{\det{F(x)}}{x}\det{F^{-1}(x)}\frac{\d
x}{2\pi i}=0.
\end{equation} However, for $m>1$ condition~(\ref{conrhp}) is not
sufficient. For example,  if $x_0\in {\cal D}$ and
\begin{equation} F(x)=\diag{x- x_0,\frac{1}{x-x_0}}\in GL(2,{\cal
O}({\mathbb C}\sm \Set{x_0})), \end{equation} then one can
 check that the corresponding Riemann-Hilbert boundary problem has
no regular solution. In order to investigate the problem further,
we shall employ  a method, based on
the Fredholm theory of linear integral equations.
It originated with D.~Hilbert
(see \cite{hilb3}, \cite{hilb2}) and was employed, in the same way as
below, by T.~Miwa in \cite{miwa}.
\begin{lmm}
\label{h1}
 Let $F(x)\in GL(m,{\cal O}({\cal A})).$ Assume that
$X^+(x),X^-(x)$ give the regular solution to the Riemann-Hilbert
boundary problem for $F(x),$ normalized by
\begin{equation}\label{norhp} X^-(\infty)=I. \end{equation} Then
$X^-(x)$ satisfies in ${\cal A}$ the linear integral equation
\begin{equation}\label{fred} X^-(x)-\oint_\gamma X^-(y)K(y,x)\frac{\d y}{2\pi
i}=I\end{equation} with kernel $K(y,x)\in gl(m,{\cal O}({\cal
A}\times{\cal A}))$ given by \begin{equation}\label{kerf}
 K(y,x)=\frac{F(y)F^{-1}(x)-I}{y-x}. \end{equation}
\end{lmm}
\begin{proof}
Using again the properties of the Cauchy integral, we obtain for
$x\in{\cal A}$ \begin{equation}  \oint_\gamma
\frac{X^-(y)F(y)-X^-(x)F(x)}{y-x}\frac{\d y}{2\pi i}= \oint_\gamma
\frac{X^+(y)-X^+(x)}{y-x}\frac{\d y}{2\pi i}=0. \end{equation}
Hence
 \begin{equation}\label{111} \oint_\gamma X^-(y)\frac{F(y)-F(x)}{y-x}\frac{\d y}{2\pi i}=
\left (\oint_\gamma\frac{X^-(x)-X^-(y)}{y-x}\frac{\d y}{2\pi
i}\right )F(x). \end{equation}
 Since
\begin{equation} \oint_\gamma\frac{X^-(x)-X^-(y)}{y-x}\frac{\d y}{2\pi
i}=X^-(x)-I, \end{equation} we obtain equation (\ref{fred}) with
kernel $K,$ given by (\ref{kerf}).
\end{proof}
Let us investigate equation~(\ref{fred}).
  According to
Fredholm theory, we should consider  Fredholm determinant
$f(\lambda)\in{\cal O}({\mathbb C}),$ defined by
\begin{equation}\label{fdet} \begin{array}{l} \displaystyle f(\lambda)=\\
\displaystyle\sum_{l=0}^\infty \frac{(-\lambda)^l}{l !}
\sum_{\alpha_1,\ldots,\alpha_l=1}^m\oint_\gamma\cdots\oint_\gamma
K\left ({x_1,\ldots, x_l\atop x_1,\ldots, x_l}\left |
{\alpha_1,\ldots,\alpha_l \atop \alpha_1,\ldots,
\alpha_l}\right.\right ) \frac{\d x_1}{2\pi i}\cdots \frac{\d
x_l}{2\pi i}, \end{array}\end{equation} where
\begin{equation}\label{str1} K\left ({x_1,\ldots, x_l\atop
y_1,\ldots, y_l}\left | {\alpha_1,\ldots, \alpha_l \atop
\beta_1,\ldots,\beta_l}\right.\right )= \det{\left
(K(x_j,y_k)_{\alpha_j\beta_k}\right )_{j,k=1}^l}. \end{equation}
If \begin{equation}\label{suf} f(1)\not=0, \end{equation} then
equation~(\ref{fred}) has unique solution $X^-(x)\in gl(m,{\cal
O}({\cal A})),$ given by \begin{equation}\label{x-}
X^-(x)=I+\frac{1}{f(1)}\oint_\gamma\hat K(y,x,1)\frac{\d y}{2\pi
i}, \quad x\in{\cal A}, \end{equation} where $\hat
K(y,x,\lambda)\in gl(m,{\cal O}({\cal A}\times{\cal
A}\times{\mathbb C}))$ is of the form \begin{equation}\label{str2}
\begin{array}{l} \displaystyle\hat K(y,x;\lambda)_{\alpha,\beta} = \\
\displaystyle\sum_{l=0}^\infty \frac{(-\lambda)^l}{l !}
\sum_{\alpha_1,\ldots,\alpha_l=1}^m\oint_\gamma\cdots\oint_\gamma
K\left ({y,x_1,\ldots, x_l\atop x,x_1,\ldots, x_l}\left |
{\alpha\alpha_1,\ldots, \alpha_l \atop
\beta,\alpha_1,\ldots,\alpha_l}\right.\right ) \frac{\d x_1}{2\pi
i}\cdots \frac{\d x_l}{2\pi i}. \end{array} \end{equation} Now we
shall show that conditions (\ref{conrhp}) and (\ref{suf}) are
sufficient in order for the Riemann-Hilbert boundary problem to
have a regular solution.
\begin{lmm}
\label{h2} Let $F(x)\in GL(m,{\cal O}({\cal A})).$ Let $K(y,x)\in
gl(m,{\cal O}({\cal A}\times{\cal A}))$ be given by (\ref{kerf})
and let $f(\lambda)\in{\cal O}({\mathbb C})$ be given by
(\ref{fdet}), (\ref{str1}). Assume that conditions~(\ref{conrhp}),
(\ref{suf}) are satisfied. Then unique solution $X^-(x)\in
gl(m,{\cal O}({\cal A}))$ of equation (\ref{fred}), given by
(\ref{x-}), (\ref{str2}), can be analytically continued in ${\cal
D}^-\cup\Set{\infty},$ is invertible there and satisfies at
$\infty$ relation (\ref{norhp}). Also, matrix function  $X^+(x)\in
gl(m,{\cal A}),$ defined by \begin{equation}\label{x+}
X^+(x)=X^-(x)F(x),\quad x\in{\cal A}, \end{equation} can be
analytically continued in ${\cal D}^+$ and is invertible there.
\end{lmm}
\begin{proof}
Since \begin{equation} \trace{K(x,x)}=
\trace{\der{F(x)}{x}F^{-1}(x)}=\der{\det{F(x)}}{x}
\det{F^{-1}(x)}, \end{equation}
 relation (\ref{conrhp})  implies \begin{equation}
\oint_{\gamma}\trace{K(y,y)}\frac{\d y}{2\pi i} =0. \end{equation}
On the other hand, it follows from (\ref{kerf}) that for
$j=0,1,\ldots$ \begin{equation}\begin{array}{l}\displaystyle
\oint_\gamma\cdots\oint_\gamma\trace{K(x_1,x_2)K(x_2,x_3) \cdots
K(x_{2j+1},x_1)} \frac{\d x_1}{2\pi i}\cdots\frac{\d
x_{2j+1}}{2\pi i}= \\\displaystyle
\oint_{\gamma}\trace{K(x,x)}\frac{\d x}{2\pi i}. \end{array}
\end{equation} In view of (\ref{str1}) we can conclude  that for
$j=0,1,\ldots$ \begin{equation} \oint_\gamma\cdots\oint_\gamma
K\left ({x_1,\ldots, x_{2j+1}\atop x_1,\ldots, x_{2j+1}}\left |
{\alpha_1,\ldots ,\alpha_{2j+1} \atop \alpha_1,\ldots,
\alpha_{2j+1}}\right.\right )\frac{\d x_1}{2\pi i}\cdots \frac{\d
x_{2j+1}}{2\pi i}=0. \end{equation} Therefore, $f(\lambda)$ is an
even function: \begin{equation} f(\lambda)=f(-\lambda).
\end{equation} Thus relation (\ref{suf}) implies \begin{equation}\label{suf-}
f(-1)\not=0. \end{equation} Let us define matrix function $Y(x)$
by \begin{equation} Y(x)=X^-(x)-I+\oint_\gamma
\frac{X^-(y)-X^-(x)}{y-x}\frac{\d y}{2\pi i}. \end{equation} Then
$Y(x)$ satisfies the linear integral system
\begin{equation}\label{fred0} Y(x)+\oint_\gamma Y(y)K(y,x)\frac{\d
y}{2\pi i}=0,\quad x\in{\cal A}. \end{equation} In view of
(\ref{suf-}) we conclude that $Y=0,$ \ie
\begin{equation}\label{x-2}
X^-(x)=I-\oint_{\gamma}\frac{X^-(y)-X^-(x)}{y-x} \frac{\d y}{2\pi
i}, \quad x\in{\cal A}. \end{equation} It follows that  we can
analytically continue $X^-(x)$ in $\ov{\mathbb C}\sm\ov{\cal D}$
with \begin{equation} X^-(x)=I-\oint_{\gamma}\frac{X^-(y)}{y-x}
\frac{\d y}{2\pi i}, \quad x\in \ov{\mathbb C}\sm\ov{\cal D}
\end{equation} and obtain
 \begin{eqnarray}
X^-(x) &\in& gl(m,{\cal O}({\cal D}^-\cup\Set{\infty})), \\
X^-(\infty)&=&I.
\end{eqnarray}
Having defined $X^+(x)\in gl(m,{\cal A})$ by (\ref{x+}), we can
conclude from (\ref{fred}) and (\ref{x-2}) that \begin{equation}
 \oint_{\gamma}\frac{X^+(y)-X^+(x)}{y-x} \frac{\d y}{2\pi i}=0, \quad x\in{\cal A}.
\end{equation} Hence we can analytically continue $X^+(x)$ in ${\cal D}$ with
 \begin{equation}
X^+(x)=\oint_{\gamma}\frac{X^+(y)}{y-x} \frac{\d y}{2\pi i},\quad
x\in{\cal D} \end{equation} and obtain \begin{equation} X^+(x)\in
gl(m,{\cal O}({\cal D}^+)). \end{equation} It remains to establish
the invertibility of $X^+,X^-.$ Relation~(\ref{conrhp}) implies
that the Riemann-Hilbert boundary problem for $\det{F(x)}$ has
regular solution $g^+(x),g^-(x),$ normalized by $g^-(\infty)=1.$
Then, according to Liouville theorem,
 \begin{eqnarray}
\det{X^+(x)}&=&g^+(x),\quad x\in {\cal D}^+, \\
\det{X^-(x)}&=&g^-(x),\quad x\in {\cal D}^-\cup\Set{\infty}.
\end{eqnarray}
Hence $X^+(x),X^-(x)$ are invertible in ${\cal D}^+,{\cal
D}^-\cup\Set{\infty},$ respectively.
\end{proof}
\begin{rem} In the
discussion above we have not made clear, whether ${\cal D}$ is a
single simply connected domain or a disjoint (with closures) union
of such. In fact, it does not matter, as long as  it is understood
that \begin{equation}
\oint_{\cup\gamma_j}=\sum_j\oint_{\gamma_j}.
\end{equation}
\end{rem}
Now let us consider the Riemann-Hilbert boundary problem, {\em
depending on a parameter.}  Let ${\cal U}\subset{\mathbb C}^n$ be
a polydisk and let $F(x,t)\in GL(m,{\cal O}({\cal A}\times{\cal
U})).$ We consider the Riemann-Hilbert boundary problem for each
$t\in{\cal U}$ and assume that for a fixed $t^0\in{\cal U}$ there
exists a regular solution, $X^+(x,t^0), X^-(x,t^0).$ In fact,
since we can always consider \begin{equation}
F^{\prime}(x,t)=X^-(x,t^0)F(x,t)(X^+(x,t))^{-1},\end{equation}{}
there is no loss of generality in the assumption
\begin{equation}\label{assrhp} F(x,t^0)\equiv I. \end{equation}
Let us observe that function
 \begin{equation}
\mu(t)=\oint_\gamma\dder{\det{F(y,t)}}{y}\det{F^{-1}(y,t)}\frac{\d
y}{2\pi i} \end{equation}
 is integer-valued and continuous,
hence constant: \begin{equation}\label{conrhpt}
\mu(t)\equiv\mu(t^0)=0.
\end{equation} Thus condition~(\ref{conrhp}) is satisfied for
every $t.$ Furthermore, kernel $K$ is holomorphic with respect to
$t$ and, because of the absolute convergence of  the Fredholm
series, the Fredholm determinant and the resolvent kernel are
holomorphic with respect to  $t:$
\begin{eqnarray} K(y,x,t)&\in& gl(m,{\cal O}({\cal A}\times{\cal A}\times{\cal U})), \\
f(\lambda,t)&\in&{\cal O}({\mathbb C}\times{\cal U}), \\
\hat K(y,x,\lambda,t)&\in&gl(m,{\cal O}({\cal A}\times{\cal
A}\times{\mathbb C}\times{\cal U})).
\end{eqnarray}
Furthermore,
assumption (\ref{assrhp}) implies that
\begin{eqnarray}
K(y,x,t^0)&\equiv& 0,\\
f(\lambda,t^0)\equiv 1.
\end{eqnarray}
Therefore, there exists  neighborhood ${\cal V}$ of $t^0,$ where
$f(1,t)$ does not vanish. It follows, in view of Lemma~\ref{h2},
 that for $t\in{\cal V}$ the Riemann-Hilbert boundary problem for $F(x,t)$
has   regular solution $X^+,X^-,$ normalized by
\begin{equation}
X^-
(\infty,t)\equiv I
\end{equation}
and
 holomorphic with respect to  $t.$
 Furthermore,
we observe that
 \begin{eqnarray}
f(1,t)X^+(x,t)&\in&gl(m,{\cal O}({\cal D}^+\times{\cal U})), \\
f(1,t)X^-(x,t)&\in&gl(m,{\cal O}(\Set{{\cal
D}^-\cup\Set{\infty}}\times{\cal U})).
 \end{eqnarray}
Since (\ref{conrhpt}) implies
\begin{eqnarray}
\det{X^+(x,t)}&\in& {\cal O}^*({\cal D}^+\times{\cal U}), \\
\det{X^-(x,t)}&\in& {\cal O}^*(\Set{{\cal
D}^-\cup\Set{\infty}}\times{\cal U}),
\end{eqnarray}
we can also observe that
\begin{eqnarray}
(f(1,t))^{m-1}(X^+(x,t))^{-1}&\in& gl(m,{\cal O}({\cal D}^+\times{\cal U})), \\
(f(1,t))^{m-1}(X^-(x,t))^{-1}&\in& gl(m,{\cal O}(\Set{{\cal
D}^-\cup\Set{\infty}}\times{\cal U})).
\end{eqnarray}
Thus we obtain
\begin{lmm}
\label{h3} Let $t^0\in{\cal U}$ and assume that matrix function
$F(x,t)\in GL(m,{\cal O}({\cal A}\times{\cal U}))$ satisfies
\begin{equation} F(x,t^0)\equiv I. \end{equation} Then there exist
neighborhood ${\cal V}$ of $t^0$ and matrix functions
\begin{eqnarray}
\label{h31}
X^+(x,t)&\in&GL(m,{\cal O}({\cal D}^+\times{\cal V})), \\
\label{h32} X^-(x,t)&\in&GL(m,{\cal O}(\Set{{\cal
D}^-\cup\Set{\infty}}\times{\cal V})),
\end{eqnarray}
such that
\begin{eqnarray}
\label{h33}
 X^-(x,t)F(x,t)&=&X^+(x,t),\quad x\in{\cal A},t\in{\cal V}, \\
\label{h34}
X^-(\infty,t)&\equiv&I.
\end{eqnarray}
Moreover, matrix functions $X^+(x,t),(X^+(x,t))^{-1},X^-(x,t),
(X^-(x,t))^{-1}$ are meromorphic
 with respect to  $t$ in ${\cal U}.$
\end{lmm}

We would like to mention that Lemma~\ref{h3} can also be proved in
another way. The existence of  neighborhood ${\cal V}$ of $t^0$
and matrix functions $X^+(x,t),X^-(x,t),$ satisfying (\ref{h31}) -
(\ref{h34}) follows immediately from H.~Grauert's theorems on
semi-continuity and on direct images of coherent sheaves under
proper analytic maps, proved in \cite{grau2}. Alternatively, it
can  be verified by utilizing  Taylor expansions with respect to
parameter $t,$ recurrence relations for Taylor coefficients and
Cauchy majorization principle. Furthermore, let
$t^1\in\partial{\cal V}$ be a fixed value of the parameter, for
which the Riemann-Hilbert boundary problem has no regular
solutions. Then we can use either A.~Grothendieck's or
G.~D.~Birkhoff's improvement of Theorem~\ref{Birthrm}, mentioned
in the previous section of this Chapter. It follows from these
results that matrix functions $X^+(x,t),(X^+(x,t))^{-1},X^-(x,t),
(X^-(x,t))^{-1}$ are meromorphic with respect to  $t$ in a
neighborhood of $t^1.$ Hence we can conclude that matrix functions
$X^+(x,t),(X^+(x,t))^{-1},X^-(x,t), (X^-(x,t))^{-1}$ are
meromorphic with respect to  $t$ in ${\cal U}.$
\chapter{The Generalized Schlesinger System}
\label{Schlesinger}
\section{Basic Notions}
We consider the linear differential system with rational
coefficients \begin{equation}\label{linsys} \der{Y}{x}=\left
(\sum_{j=1}^n\sum_{k=0}^{p_{j}} \frac{Q_{j,k}}{(x-
t_{j})^{k+1}}\right )Y, \end{equation} where $Q_{j,k}\in
gl(m,{\mathbb C}),$ and $t_1,\ldots,t_n$ are fixed mutually
distinct points in ${\mathbb C}.$ To simplify our presentation, we
assume\footnote{There is no loss of generality in this assumption.
One can  use a M\"obius transformation to ensure that it is
satisfied, although this will increase $n$ by $1.$}  that
system~(\ref{linsys}) is regular at $x=\infty,$ \ie that
\begin{equation}\label{reg} \sum_{j=1}^n Q_{j,0}=0. \end{equation}
According to Cauchy theorem, in a neighborhood of $x=\infty$ there
exists  unique holomorphic fundamental solution $Y$  of
system~(\ref{linsys}), satisfying the initial condition
\begin{equation}\label{inisys} Y(\infty)=I.
\end{equation} Moreover, $Y(x)$ admits the analytic continuation
along any path in \begin{equation} {\cal R}_t=\ov{\mathbb
C}\sm\Set{ t_{j}}_{j=1}^n, \end{equation} staying invertible.
Since for $n>1$ ${\cal R}_t$ is {\em not} simply connected, the
result of the analytic continuation depends, in general, on the
homotopy class of the path along which it is performed. Therefore,
we should consider the universal covering \begin{equation}
\psi_t:\tilde{\cal R}_t\mapsto{\cal R}_t. \end{equation} We denote
a point of $\tilde{\cal R}_t$ in $\psi^{-1}(x)$ by $\tilde x$ and
the group of deck transformations, acting on $\tilde{\cal R}_t,$
by $\Delta_t.$ We recall that $\tilde{\cal R}_t$ is a simply
connected Riemann surface with the complex differential structure
provided by $\psi_t.$  If one distinguishes a point in $\psi_t^{-
1}(\infty)$ then $\tilde x$ can be  interpreted as a homotopy
class of paths in ${\cal R}_t,$ starting at $\infty$ and ending at
$x.$ This interpretation allows us to identify $\Delta_t$ with
fundamental group $\pi_1({\cal R}_t,\infty),$ the action being the
usual multiplication of homotopy classes. Namely, for
$\sigma\in\Delta_t$ and $\tilde x\in\tilde{\cal R}_t$ we interpret
$\sigma\tilde x\in\psi_t^{-1}(x)$ as the homotopy class of the
path, obtained by  traversing first a loop in class $\sigma,$
beginning and ending at $\infty,$  and then a representative of
class $\tilde x$ from $\infty$ to $x.$ Thus we can define
\begin{equation} Y(\tilde x)\in GL(m,{\cal O}(\tilde{{\cal R}_t}),
\end{equation} satisfying system~(\ref{linsys}) in $\tilde{\cal
R}_t,$ and initial condition~(\ref{inisys}), where $\infty$ now
denotes the distinguished point. Let $\sigma\in\Delta_t$ then
$Y\circ\sigma(\tilde x)$ is another fundamental solution of
system~(\ref{linsys}) in $\tilde{\cal R}_t,$ determined by the
initial condition \begin{equation}
Y\circ\sigma(\infty)=Y(\sigma\infty).
\end{equation} Hence we can consider mapping \begin{equation}\label{monodromy} \Phi:
\Delta_t\mapsto GL(m,{\mathbb C}), \end{equation} defined by
\begin{equation}\label{defmon} \Phi(\sigma)=Y^{-1}(\sigma\infty)\in GL(m,{\mathbb C})
\end{equation} and satisfying \begin{equation}\label{defmon2} ( Y\circ\sigma(\tilde
x))\Phi(\sigma)\equiv Y(\tilde x). \end{equation} Note that for
any $\sigma^{\prime},\sigma^{\prime\prime}\in \Delta_t,$
\begin{equation}\begin{array}{l}
Y\circ(\sigma^{\prime}\sigma^{\prime\prime})=(Y\Phi(\sigma^{\prime})^{-
1})\circ\sigma^{\prime\prime}= \\
(Y\circ\sigma^{\prime\prime})\Phi(\sigma^{\prime})^{-
1}=Y\Phi(\sigma^{\prime\prime})^{-1} \Phi(\sigma^{\prime})^{-1},
\end{array} \end{equation} \ie \begin{equation}
\Phi(\sigma^{\prime})\Phi(\sigma^{\prime\prime})=\Phi(\sigma^{\prime}\sigma^{\prime\prime}),
\end{equation} and \begin{equation} \Phi(1)=I. \end{equation} We conclude that $\Phi$ is a linear
representation of group $\Delta_t.$ \begin{dfn} We shall call
representation $\Phi$ the  monodromy representation of $Y.$
\end{dfn}  Now we would like to introduce a certain canonical
factorization of $Y(\tilde x)$ for $x$ in a simply connected open
neighborhood  of  singular point $t_j.$ In order to do this, we
need to introduce some notations first.

 For $j=1,\ldots,n$ let
${\cal L}_j$ be a simply connected open neighborhood  of $t_j,$
such that \begin{equation} {\cal L}_{j^{\prime}}\cap
t_{j^{\prime\prime}}=\emptyset,\quad
j^{\prime}\not=j^{\prime\prime}. \end{equation} Denote
\begin{equation} {\cal L}_{t_j}={\cal L}_j\sm t_j\subset{\cal
R}_t. \end{equation} Also, let $x_j\in{\cal L}_{t_j},\tilde
x_j\in\tilde{\cal R}_t,$ and denote by $\tilde{\cal L}_{t_j}$  the
connected component\footnote{For $n>2$ set $\psi_t^{-1}({\cal L}_{
t_{j}})$ is disconnected.}  of $\psi_t^{- 1}({\cal L}_{ t_{j}}),$
containing $\tilde x_j.$ Then $\tilde x\in\tilde{\cal L}_{t_j}$
if, and only if, the homotopy class $\tilde x_j^{-1}\tilde x$ of a
path in ${\cal R}_t,$ starting at $x_j$ and ending at $x,$  has a
representative, which lies entirely in ${\cal L}_{t_j}.$ Moreover,
such a representative is unique up to homotopy in ${\cal
L}_{t_j}.$ It follows that $\tilde{\cal L}_{t_j}$ is simply
connected and that \begin{equation}\label{loccov}
\psi_t\val{\tilde{\cal L}_{t_j}}:\tilde{\cal L}_{t_j}\mapsto {\cal
L}_{t_j}
\end{equation} is the universal covering over ${\cal L}_{t_j}.$
The group of deck transformations of covering~(\ref{loccov}),
which we denote by $\Delta_{t_j},$ is a cyclic subgroup of
$\Delta_t.$ Its generator $\sigma_j$ is defined in $\pi_1({\cal
R}_t,\infty)$  by \begin{equation}\label{defsi} \sigma_j=\tilde
x_j\varsigma_j\tilde x_j^{-1}, \end{equation} where
$\varsigma_j\in\pi_1({\cal L}_{t_j},x_j)$ is the homotopy class of
a simple loop $\gamma_j,$ making one positive circuit of $t_j.$
Fundamental group $\pi_1({\cal L}_{t_j},x_j)$ can be identified
with $\pi_1({\cal R}_{t_j},x_j),$ where
\begin{equation}\label{rtj} {\cal R}_{t_j}={\mathbb C}\sm t_j.
\end{equation} Denote by \begin{equation}\label{trtj} \psi_{ t_{j}}:\tilde{\cal R}_{
t_{j}}\mapsto{\cal R}_{ t_{j}} \end{equation}  the universal
covering over ${\cal R}_{ t_{j}},$ then $\tilde{\cal L}_{ t_{j}}$
can be identified with $\psi_{ t_{j}}^{- 1}({\cal L}_{
t_{j}})\subset\tilde{\cal R}_{ t_{j}},$ and $\Delta_{ t_{j}}$ -
with the group of deck transformations of covering~(\ref{trtj}).
Furthermore, we can assume that for $j=1,\ldots,n$ homotopy class
$\tilde x_j$ has a representative, which is simple and has no
points of intersection with $\gamma_j,$ except $x_j.$ Then
$\sigma_1,\ldots,\sigma_n,$ defined by (\ref{defsi}),
generate\footnote{ These are not free generators; taking more care
in choice of $\tilde x_j,$ we can always obtain, for example,
$\sigma_1\cdots\sigma_n=1.$}  group $\pi_1({\cal R}_t,\infty).$
Hence subgroups $\Delta_{t_1},\ldots\Delta_{t_n}$ generate group
$\Delta_t,$ and the {\em monodromy generators}
$\Phi(\sigma_1),\ldots,\Phi(\sigma_n)$ completely determine the
monodromy representation of $Y.$  Let $\tilde{{\mathbb C}^*}$
denote the Riemann surface of the logarithm, \ie the universal
covering surface over ${\mathbb C}^*={\mathbb C}\sm0.$ Considering
$\tilde{\cal R}_{t_j},$ with the projection \begin{equation}
\psi_{t_j}-t_j:\tilde{\cal R}_{t_j}\mapsto{\mathbb C}^*,
\end{equation} as a covering surface over ${\mathbb C}^*,$ we can
observe that there exists an isomorphism of covering spaces,
mapping $\tilde{\cal R}_{ t_{j}}$ onto $\tilde{\mathbb C}^*.$ Such
an isomorphism is not unique, but we shall fix for the moment some
arbitrary one and denote it by \begin{equation} \tilde
x\mapsto\tilde x- t_{j}. \end{equation} This allows us to consider
\begin{equation} \log(\tilde x-t_j)\in{\cal O}(\tilde{\cal R}_{t_j}).
\end{equation} Let $\sigma_0$ be the generator of the group of
deck transformations, acting on $\tilde{\mathbb C}^*$ and
corresponding to one positive circuit of $0.$ Then
\begin{equation} \log\sigma_0\tilde x=\log\tilde x+2\pi i,\quad
\tilde x\in\tilde{\mathbb C}^* \end{equation} and
\begin{equation}\label{sigma0} \sigma_0(\tilde x-
t_j)=\sigma_j\tilde x -t_j, \quad x\in\tilde{\cal R}_{t_j}.
\end{equation} Thus we have \begin{equation}\label{log} \log(\sigma_j\tilde
x-t_j)=\log(\tilde x- t_j)+2\pi i.
\end{equation}

With these notations we can  formulate the following result,
due to G.~D.~Birkhoff (see \cite{bir1}):
\begin{thrm}[Birkhoff] \label{faclmm} Let $Y(\tilde x)\in
GL(m,{\cal O}(\tilde{\cal R}_t))$ be the fundamental solution of
system~(\ref{linsys}) with initial condition (\ref{inisys}). Then
for $j=1,\ldots,n$ $Y(\tilde x)$  admits in $\tilde{\cal L}_{
t_{j}}$ the  factorization \begin{equation}\label{fac} Y(\tilde
x)=H_{j}(x)M_{j}(\tilde x- t_{j}), \end{equation} where
\begin{eqnarray}
H_{j}(x) &\in & GL(m,{\cal O}({\cal L}_{j})), \label{anh} \\
M_{j}(\tilde x) & \in & GL(m,{\cal O}(\tilde{\mathbb C}^*)).
\label{anm}
\end{eqnarray}
\end{thrm}
\begin{proof} Let us fix $j.$ Since monodromy generator
$\Phi(\sigma_j)$ is non-degenerate, we can choose some value of
the matricial logarithm $\log\Phi(\sigma_j)$ and consider
\begin{equation} X_{j}(\tilde x)=\tilde x^{\frac{1}{2\pi
i}\log\Phi(\sigma_j)}= \exp(\frac{1}{2\pi
i}\log\Phi(\sigma_j)\log\tilde x) \in GL(m,{\cal O}(\tilde{C^*})).
\end{equation} In view of (\ref{log}), $X_{j}(\tilde x- t_{j})\in
GL(m,{\cal O}(\tilde{\cal R}_{t_j}))$ satisfies \begin{equation}
X_{j}(\sigma_j\tilde x- t_{j})=\Phi(\sigma_j)X_{j}(\tilde x-
t_{j}). \end{equation} Hence $Y(\tilde x)X_{j}(\tilde x- t_{j})$
is actually univalued in ${\cal L}_{ t_{j}},$ and we can define
\begin{equation} F(x)= \left (Y(\tilde x)X_{j}(\tilde x- t_{j})\right
)^{-1}\in GL(m,{\cal O}({\cal L}_{t_j})). \end{equation} Now we
can consider the Riemann-Hilbert boundary problem for $F(x),$
which we have discussed in Chapter \ref{intro}. By
Theorem~\ref{Birthrm}, there exist
\begin{eqnarray} X^+(x)&\in&GL(m,{\cal O}({\cal L}_j)), \\
X^-(x)&\in&GL(m,{\cal O}({\cal R}_{t_j})),
\end{eqnarray} such that \begin{equation} X^-(x)F(x)=X^+(x),\quad
x\in{\cal L}_{t_j}. \end{equation} It suffices to define
\begin{eqnarray}
H_j(x)&=&(X^+)^{-1}(x)\in GL(m,{\cal O}({\cal L}_j)), \\
M_{j}(\tilde x)&=&X^-(x+ t_{j})X_{j}^{-1}(\tilde x)\in GL(m,{\cal
O}(\tilde{\mathbb C}^*))
\end{eqnarray} in order to complete the proof. \end{proof}
\begin{dfn} We shall call $H_{j}(x)$ and $M_{j}(\tilde x),$
respectively, the  non-singular and principal factors of $Y$ at
$t_j.$ \end{dfn} \begin{rem} \label{facrem} \begin{enumerate}
\item Of course, factors $H_j$ and $M_j$ are defined up to the
transformation \begin{equation} \left\{\begin{array}{l}H_j^{\prime}(x)=H_j(x)E_j(x-t_j),\\
 M_j^{\prime}(\tilde x)=E_j^{-1}(x)M_j(\tilde x),\end{array}\right. \end{equation}
 for arbitrary $E_j(x)\in
GL(m,{\cal O}({\mathbb C})).$ \item Although our choice of
$\tilde{\cal L}_{t_j}$ and $\tilde x-t_j$ is somewhat ambiguous,
one should keep in mind that any other connected component of
$\psi_t^{-1}({\cal L}_{t_j})$ is of the form $\sigma\tilde{\cal
L}_{t_j}$ for $\sigma\in\Delta_t,$ and any other isomorphism,
mapping $\tilde{\cal R}_{t_j}$ onto $\tilde{\mathbb C}^*,$ is of
the form $\tilde x\mapsto\sigma_j^l\tilde x- t_j$ for
$l\in{\mathbb Z}.$ According to (\ref{defmon2}), relation
(\ref{fac}) implies \begin{equation}\label{prinmon}
M_j(\sigma_j\tilde x- t_j)=M_j(\sigma_0(\tilde x-t_j))= M_j(\tilde
x- t_j)\left (\Phi(\sigma_j)\right )^{-1}. \end{equation} Hence we
can fix non-singular factor $H_j(x)$  independently of these
choices and define for arbitrary $\sigma\in\Delta_t$ the
corresponding principal factor by \begin{equation}
M_j(\sigma\tilde x-t_j)=M_j(\tilde x- t_j)\left
(\Phi(\sigma)\right )^{-1},
\end{equation} so that \begin{equation} H_j(x)M_j(\sigma\tilde
x-t_j)=Y(\sigma\tilde x)=Y(\tilde x)\Phi(\sigma^{-1}).
\end{equation} holds true. \end{enumerate}
\end{rem} We describe some basic properties of the non-singular
and principal factors in the following \begin{lmm}
\label{propfaclmm}
\begin{enumerate} \item Let $Y(\tilde x)\in GL(m,{\cal O}(\tilde{\cal R}_t))$
be the fundamental solution of system~(\ref{linsys}) with initial
condition (\ref{inisys}).  Let $H_j,M_j$ be the non-singular and
principal factors of $Y$ at $t_j.$ Then  there exists $P_j(x)\in
gl(m,{\cal O}({\mathbb C})),$ such that $H_j$ and $M_j$ satisfy
the systems
\begin{equation}\label{holsys}\begin{array}{rl}
\displaystyle\der{H_{j}}{x}=& \displaystyle\left
(\sum_{j^{\prime}=1}^n\sum_{k=0}^{p_{j^{\prime}}}
\frac{Q_{j^{\prime},k}}{(x- t_{j^{\prime}})^{k+1}}\right )H_{j}(x)- \\
&\displaystyle H_j(x) \left (\sum_{k=0}^{p_{j}}\frac{
J_{j,k}}{(x-t_j)^{k+1}}+P_j(x-t_j)\right ), \quad x\in{\cal L}_{
t_{j}}
\end{array} \end{equation} and \begin{equation}\label{linm}
\der{M_j}{x}=\left (\sum_{k=0}^{p_{j}}\frac{
J_{j,k}}{x^{k+1}}+P_j(x)\right )M_j(\tilde x), \end{equation}
where $J_{j,k}\in gl(m,{\mathbb C})$ is given by
\begin{equation}\label{atoj} J_{j,k} =\sum_{k^{\prime}=0}^{p_{j}-
k}\oint_{\gamma_j}\frac{ H_{j}^{-1}(x)
Q_{j,k+k^{\prime}}H_{j}(x)}{(x- t_j)^{k^{\prime}+1}}\frac{\d
x}{2\pi i}. \end{equation} \item Let $Y(\tilde x)\in GL(m,{\cal
O}(\tilde{\cal R}_t))$ be the fundamental solution of
system~(\ref{linsys}) with initial condition (\ref{inisys}). Let
$J_{j,k}\in gl(m,{\mathbb C}),$ $P_j(x)\in gl(m,{\cal O}({\mathbb
C}))$ be such that system~(\ref{holsys}) has solution $H_j(x)\in
GL(m,{\cal O}({\cal L}_j)).$ Then $H_j$ is the non-singular factor
of $Y$ at $t_j,$ and corresponding principal factor $M_j$ is a
fundamental solution of system~(\ref{linm}). \item Assume that
$Y(\tilde x)\in GL(m,{\cal O}(\tilde{\cal R}_t))$ satisfies
(\ref{inisys}) and, for $j=1,\ldots,n$  admits in $\tilde{\cal
L}_{ t_{j}}$ factorization (\ref{fac}), where $H_j(x)\in
GL(m,{\cal O}({\cal L}_j))$ and $M_j(\tilde x)$ is a fundamental
solution of system~(\ref{linm}) with some $J_{j,k}\in
gl(m,{\mathbb C}),$ $P_j(x)\in gl(m,{\cal O}({\mathbb C})).$ Then
$Y$ is a fundamental solution of system (\ref{linsys}), where
$Q_{j,k}\in gl(m,{\mathbb C})$ are given by
\begin{equation}\label{jtoa} Q_{j,k} =\sum_{k^{\prime}=0}^{p_{j}-
k}\oint_{\gamma_j}\frac{ H_j(x) J_{j,k+k^{\prime}}H_j^{-1}(x)}{(x-
t_j)^{k^{\prime}+1}}\frac{\d x}{2\pi i}. \end{equation}
\end{enumerate}
\end{lmm}
\begin{proof} The proof is done by straightforward computation.
We shall prove statement 1, and the rest can be done absolutely
analogously.

Let $Y(\tilde x)\in GL(m,{\cal O}(\tilde{\cal R}_t))$ be the
fundamental solution of system~(\ref{linsys}) with initial
condition (\ref{inisys}) and let $H_j,M_j$ be the non- singular
and principal factors of $Y$ at $t_j.$ Then it follows from
(\ref{prinmon}) that  the logarithmic derivative of $M_j(\tilde
x),$ \begin{equation}  \der{M_j}{x}M_j^{-1}(\tilde x),
\end{equation}  is univalued in $C^*.$ In view of (\ref{linsys})
and (\ref{fac}), it has at $x=0$ a pole of order $p_j+1.$ Thus
$M_j(\tilde x)$ is a fundamental solution of system~(\ref{linm})
where $J_{j,k}\in gl(m,{\mathbb C}),$ $P_j(x)\in gl(m,{\cal
O}({\mathbb C})).$ Substituting (\ref{linm}) into (\ref{linsys}),
(\ref{fac}), we conclude that $H_j(x)$ satisfies system
(\ref{holsys}) and that $J_{j,k}$ are given by (\ref{atoj}).
\end{proof} \begin{rem}
\begin{enumerate} \item If we retrace our proof of Theorem
\ref{faclmm} and use there Theorem \ref{Birthrm}  to  demand from
$X^-$ to have at most a pole at $\infty$, then we can conclude
that also the logarithmic derivative of $M_j$ has at most a pole
at $\infty.$ (Actually, this stronger formulation of
Theorem~\ref{faclmm} appears in G.~D.~Birkhoff`s work
\cite{bir1}.) Hence $P_j(x)$ in system~(\ref{linm}) can always be
chosen to be a polynomial. \item We would like to note here that
one should be able, in principle, to recover  $Q_{j,k}$ from
principal factors $M_j$ alone. Indeed, it follows from Liouville
theorem that if for $j=1,\ldots,n$ $M_j\in GL(m,{\cal
O}(\tilde{\mathbb C}^*))$ satisfy monodromy relations
(\ref{prinmon}) then there exists at most one  $Y \in GL(m,{\cal
O}(\tilde{\cal R}_t)),$ normalized by (\ref{inisys}), which admits
in $\tilde{\cal L}_{t_j}$ factorizations (\ref{fac}) with some
$H_j(x)\in GL(m,{\cal O}({\cal L}_j)).$ \end{enumerate} \end{rem}
\section{The Generalized Schlesinger
System} Let us consider an analytic family of linear systems  of
type~(\ref{linsys}), parameterized by position of the singular
points.  Let us denote \begin{equation} {\cal S}={\mathbb
C}^n\sm\bigcup_{j^{\prime}\not=j^{\prime\prime}}\Set{t=(t_1,\ldots,t_n):t_{j^{\prime}
} =t_{j^{\prime\prime}}} \end{equation} and let ${\cal U}={\cal
U}_1\times\ldots\times{\cal U}_n\subset{\mathbb C}^n$ be a
polydisk, such that $\ov{\cal U}\subset{\cal S}.$ Assume that
$Q_{j,k}(t)\in gl(m,{\cal O}({\cal U}))$ satisfy
\begin{equation}\label{regt} \sum_{j=1}^{n} Q_{j,0}(t)\equiv 0
\end{equation} and consider for each fixed $t=(t_1,\ldots,t_n)
\in{\cal U}$ the system \begin{equation}\label{linsyst}
\der{Y}{x}=\sum_{j=1}^n\sum_{k=0}^{p_{j}} \frac{Q_{j,k}(t)}{(x-
t_{j})^{k+1}}Y \end{equation} with the initial condition
\begin{equation}\label{inidef} Y(\infty,t)\equiv I. \end{equation} We encounter here
a slight difficulty, because for different $t$ we have to consider
$\tilde x$ on different surfaces $\tilde {\cal R}_t.$ However, it
can be overcome if we consider \begin{eqnarray}
{\cal T}_{j}&=&\Set{(x,t)\in\ov{\mathbb C}\times{\cal U}:x= t_{j}}, \\
{\cal R}&=&\ov{\mathbb C}\times{\cal U}\sm\bigcup_{j}{\cal T}_{j},
\end{eqnarray} and the universal covering \begin{equation} \psi:\tilde{\cal
R}\mapsto{\cal R}, \end{equation} with the group of deck
transformations $\Delta.$ Our reasoning here is analogous to what
we did in the previous section, so we shall be as brief as
possible.

 First of all, for a fixed $t\in{\cal U}$  let us
denote \begin{eqnarray}
({\cal R}_{t},t)&=&\Set{(x,t):x\in{\cal R}_{t}}\subset{\cal R}, \\
(\tilde{\cal R}_{t},t)&=&\psi^{-1}({\cal R}_{t},t). \end{eqnarray}
Then set $(\tilde{\cal R}_{t},t)$ is simply  connected, hence
isomorphic to $\tilde{\cal R}_{t},$ and  we can, indeed, interpret
a point of $\psi^{-1}(x,t)$ as a pair $(\tilde x,t),$ where
$\tilde x\in\tilde{\cal R}_t.$ Analogously, for a  fixed
$x\in\ov{\mathbb C}\sm\bigcup_{j=1}^n{\cal U}_j$ we denote
\begin{equation} (x,{\cal U})=\Set{(x,t):t\in{\cal U}}\subset{\cal
R}
\end{equation} and observe that the set $\psi^{-1}(x,{\cal U})$ is
disconnected\footnote{ We assume, of course, $n>1.$}. Let us
choose a connected component, and denote it by $(\tilde x,{\cal
U}).$ Then   $(\tilde x,{\cal U})$ is simply connected (and,
therefore, isomorphic to ${\cal U}$), hence for any $t\in{\cal U}$
the intersection $(\tilde x,{\cal U})\cap(\tilde{\cal R}_t,t)$
consists of a single point $(\tilde x,t).$ Thus we can identify
$\tilde x\in\tilde{\cal R}_{t^{\prime}}$ and $\tilde x\in
\tilde{\cal R}_{t^{\prime\prime}},$ such that $(\tilde
x,t^{\prime})$ and $(\tilde x,t^{\prime\prime})$ belong to the
same connected component $(\tilde x,{\cal U})$ of $\psi^{-
1}(\Set{(x,t):t\in{\cal U}}).$ In particular, if we assume that
initial condition (\ref{inidef}) is taking  place at some such
connected component of $\psi^{-1}(\infty,t)$ then, because of
analytic dependence of solution of a linear system on the
coefficients, the family of systems~(\ref{linsyst}) determines
$Y(\tilde x,t)\in GL(m,{\cal O}(\tilde{\cal R}),$  and, instead of
(\ref{linsyst}), we can actually write
\begin{equation}\label{deform}
\dder{Y}{x}=\sum_{j=1}^n\sum_{k=0}^{p_{j}} \frac{Q_{j,k}(t)}{(x-
t_{j})^{k+1}}Y. \end{equation} Also, since  group of deck
transformations $\Delta$  acts on each surface $(\tilde{\cal
R}_{t},t),$ it can be identified with the corresponding group of
deck transformations $\Delta_t,$ with the understanding
\begin{equation} \sigma(\tilde x,t)=(\sigma\tilde x,t),
\end{equation}  and thus the monodromy representations of $Y$ form
an analytic family, that is, for $\sigma\in \Delta$ we can
consider the monodromy matrix function \begin{equation}
\Phi(\sigma)(t)=\left (Y\circ\sigma(\tilde x,t)\right
)^{-1}Y(\tilde x,t)= Y^{-1}(\sigma\infty, t)\in GL(m,{\cal
O}({\cal U})).
\end{equation}
\begin{dfn} We shall say that system~(\ref{deform}) gives an
isomonodromic deformation if the monodromy representation of $Y$
is independent of $t,$ \ie for each  $\sigma\in\Delta$ the
monodromy matrix function $\Phi(\sigma)$ is constant in ${\cal
U}.$
\end{dfn} Let us consider now a polydisk
${\cal L}={\cal L}_1\times\ldots\times{\cal L}_n\subset{\cal S},$
such that $\ov{\cal U}\subset{\cal L},$ and for $j=1,\ldots,n$ let
us denote
\begin{eqnarray} {\cal L}_{{\cal T}_{j}}&=&{\cal L}_{j}\times{\cal U}\sm{\cal T}_{j}, \\
{\cal R}_{{\cal T}_{j}} &=&{\mathbb C}\times{\cal U}\sm{\cal
T}_{j}.
\end{eqnarray} As before, we can choose a connected component
$\tilde{\cal L}_{{\cal T}_{j}}$ of $\psi^{- 1}({\cal L}_{{\cal
T}_{j}})$ in such a way that the corresponding cyclic subgroups
$\Delta_{{\cal T}_{j}}$ of $\Delta,$ acting on $\tilde{\cal
L}_{{\cal T}_{j}},$ generate $\Delta,$ and fix an isomorphism
\begin{equation} (\tilde x,t)\mapsto(\tilde x- t_{j},t), \end{equation} which
maps the universal covering space $\tilde{\cal R}_{{\cal T}_{j}}$
over ${\cal R}_{{\cal T}_{j}}$ onto $\tilde{\mathbb
C}^*\times{\cal U}.$
\begin{dfn} We shall say that system~(\ref{deform}) gives the {\em
isoprincipal deformation}, if  $Y$ has the principal factors,
independent of $t,$ \ie for $j=1,\ldots,n$ $Y$ admits in
$\tilde{\cal L}_{{\cal T}_{j}}$ the factorization
\begin{equation}\label{fact} Y(\tilde x,t)=H_{j}(x,t)M_{j}(\tilde
x- t_{j}), \end{equation} where
\begin{eqnarray}
\label{anht} H_{j}(x,t)&\in& GL(m,{\cal O}({\cal L}_{j}\times{\cal U})),\\
\label{anmt} M_{j}(\tilde x)&\in& GL(m,{\cal O}(\tilde{\mathbb
C}^*)).
\end{eqnarray} \end{dfn} \begin{rem} \label{schliso} In view of
relation (\ref{prinmon}) in Remark~\ref{facrem}, the isoprincipal
deformation belongs to the class of isomonodromic deformations.
\end{rem} Now we shall try to develop  some simple criteria of
whether given system (\ref{deform}) gives the isoprincipal
deformation. \begin{lmm} \label{schlmm} System~(\ref{deform})
gives the isoprincipal deformation if, and only if, $Y$ satisfies
the  linear Pfaffian system \begin{equation}\label{isoschl} \d Y=
\sum_{j=1}^n\sum_{k=0}^{p_{j}} Q_{j,k}(t)\frac{\d x-\d t_{j}}{(x-
t_{j})^{k+1}}Y. \end{equation} \end{lmm}
\begin{proof} Let us assume first that system~(\ref{deform}) gives
the isoprincipal deformation. Since $Y(\tilde x,t)$ satisfies,  by
definition, system (\ref{deform}), we only have to prove that
\begin{equation}\label{dytj}
\dder{Y}{t_j}Y^{-1}=-\sum_{k=0}^{p_j}\frac{Q_{j,k}(t)}{(x-
t_j)^{k+1}}. \end{equation} Let us observe that, since the
isoprincipal
 deformation is isomonodromic, for any
$\sigma\in\Delta$ \begin{equation} \left
(\dder{Y}{t_j}Y^{-1}\right )\circ\sigma=
\dder{Y}{t_j}Y^{-1}-Y\left (\Phi(\sigma^{-1})\dder{\Phi(\sigma)}{
t_{j}}\right )Y^{-1} =\dder{Y}{t_j}Y^{-1}, \end{equation} hence,
\begin{equation} \dder{Y}{t_j}Y^{-1}\in gl(m,{\cal O}({\cal R})).
\end{equation} Taking into account factorization (\ref{fact}) of
$Y$ for $j^{\prime}\not=j,$ we can observe that in ${\cal
L}_{{\cal T}_{j^{\prime}}}$ \begin{equation} \dder{Y}{t_j}Y^{-
1}=\dder{H_{j^{\prime}}}{t_j}H_{j^{\prime}}^{-1} \in gl(m,{\cal
O}({\cal L}_{j^{\prime}}\times{\cal U})), \end{equation} hence,
\begin{equation} \dder{Y}{t_j}Y^{- 1}\in gl(m,{\cal O}({\cal R}_{{\cal
T}_j})). \end{equation} Analogously, since \begin{equation}
\dder{(M_{j}(\tilde x- t_{j}))}{t_j}=-\dder{(M_{j}(\tilde x-
t_{j}))}{x}, \end{equation} we can observe that in ${\cal
L}_{{\cal T}_{j}}$ \begin{equation} \left
(\dder{Y}{t_j}+\dder{Y}{x}\right )Y^{-1}= \left
(\dder{H_{j}}{t_j}+\dder{H_{j}}{x}\right ) H_j^{-1}\in gl(m,{\cal
O}({\cal L}_j\times{\cal U})),
\end{equation} hence \begin{equation} \dder{Y}{t_j}Y^{-
1}+\sum_{k=0}^{p_{j}}\frac{Q_{j,k}(t)}{(x- t_{j})^{k+1}}\in
gl(m,{\cal O}(\ov{\mathbb C}\times{\cal U})). \end{equation} But
initial condition (\ref{inidef}) implies
\begin{equation}\label{dreg} \d Y\val{\tilde x=\infty}= 0,
\end{equation} and, applying Liouville theorem, we can conclude
\begin{equation} \dder{Y}{t_j}Y^{-
1}+\sum_{k=0}^{p_{j}}\frac{Q_{j,k}(t)}{(x- t_{j})^{k+1}} \equiv 0.
\end{equation} Thus, indeed,  (\ref{dytj}) holds true and  $Y$ satisfies
system~(\ref{isoschl}).

Conversely, if $Y$ satisfies system~(\ref{isoschl}), then let us
choose a simply connected neighborhood ${\cal V}$ of $0$ in
${\mathbb C},$ such that \begin{equation} {\cal
V}\subset\Set{{\cal L}_{j}- t_{j}}\quad\fa t\in{\cal U} ,\fa j
\end{equation} and consider for any fixed $\tilde v\in\tilde{\cal
V}^*$ the surface $\Set{\tilde x=\tilde v+ t_{j}},$ on which $Y$
satisfies \begin{equation} \d
Y=\sum_{j^{\prime}\not=j}\sum_{k=0}^{p_{j^{\prime}}}
Q_{j^{\prime},k}\frac{\d t_{j}-\d t_{j^{\prime}}}{(v+
t_{j}-t_{j^{\prime}})^{k+1}}Y. \end{equation} Hence, the linear
Pfaffian system \begin{equation}\label{loccomp} \d
H_{j}=\sum_{j^{\prime}\not=j}\sum_{k=0}^{p_{j^{\prime}}}
Q_{j^{\prime},k}\frac{\d t_{j}-\d t_{j^{\prime}}}{(v+
t_{j}-t_{j^{\prime}})^{k+1}}H_{j}, \end{equation} where $v\in{\cal
V}$  is a parameter, is  compatible\footnote{Compatibility of a
Pfaffian system means that  the one-form on the right-hand side is
closed.} and thus, according to Frobenius theorem, completely
integrable. Let us fix $t^0\in{\cal U}$ and let $H_{j}(x,t^0)\in
GL(m,{\cal O}({\cal L}_{j}))$ and $M_{j}(\tilde x)\in GL(m,{\cal
O}(\tilde{\mathbb C}^*))$ be the non-singular and principal
factors of  $Y(\tilde x,t^0).$ Because of the linearity of
system~(\ref{loccomp}), it has a solution  $ H_{j}(v,t)\in
GL(m,{\cal O}({\cal V}\times{\cal U})),$ satisfying
\begin{equation}\label{varini} H_{j}(v,t^0)=H_{j}^0(v + t_{j}^0).
\end{equation} The product ${H_{j}}^{- 1}(v,t)Y(\tilde v+
t_{j},t)\in GL(m,{\cal O}(\tilde{\cal V}^*\times{\cal U}))$ is
independent of $t,$ hence \begin{equation}\label{defm}
{H_{j}}^{-1}(v,t)Y(\tilde v+ t_{j},t)=M_{j}(\tilde v)
\end{equation} and $H_{j}(v,t)$ can be analytically continued into the simply
connected domain $\Set{(v,t):t\in{\cal U},v\in{\cal L}_{j}-
t_{j}}\subset{\mathbb C}^*\times{\cal U}.$ Returning to $\tilde
x=\tilde v+ t_{j},$ we obtain that $H_{j}(x,t)\in GL(m,{\cal
O}({\cal L}_{j}\times{\cal U}))$ and \begin{equation} Y(\tilde
x,t)=H_{j}(x,t)M_{j}(\tilde x- t_{j}), \quad (\tilde
x,t)\in\tilde{\cal L}_{{\cal T}_{j}}, \end{equation} \ie
system~(\ref{deform}) gives the isoprincipal deformation, and this
completes the proof.
\end{proof} Let us observe that, since  system (\ref{isoschl})
implies (\ref{deform}) and (\ref{dreg}),  Lemma \ref{schlmm}
actually means that system~(\ref{deform}) gives the isoprincipal
deformation if, and only if, system~(\ref{isoschl}) is integrable,
\ie compatible. The compatibility condition is
\begin{equation}\label{compschl}
\sum_{j=1}^n\sum_{k=0}^{p_{j}}\frac{\d Q_{j,k}\wedge\d(x- t_{j})}
{(x- t_{j})^{k+1}}=\sum_{j,j^{\prime}=1}^n\sum_{k=0}^{p_{j}}
\sum_{k^{\prime}=0}^{p_{j^{\prime}}}\frac{Q_{j,k}Q_{j^{\prime},k^{\prime}}\d(x-
t_{j})\wedge\d(x-t_{j^{\prime}})} {(x- t_{j})^{k+1}(x-
t_{j^{\prime}})^{k^{\prime}+1}}. \end{equation} Since $Q_{j,k}$
are independent of $x,$ both sides of (\ref{compschl}) are
rational with respect to  $x.$ Comparing on both sides the
principal parts of Laurent series at each pole, we obtain the
following {\em non-linear} Pfaffian system:
\begin{equation}\label{schl}
\begin{array}{l} \d Q_{j,k}=\sum_{l=0}^{p_{j}-k}\sum_{j^{\prime}\not=j}
\sum_{k^{\prime}=0}^{p_{j^{\prime}}}(-1)^l{l+k^{\prime}\choose l}
\left  [Q_{j^{\prime},k^{\prime}},Q_{j,k+l}\right  ]\frac{\d(
t_{j}-t_{j^{\prime}})}{( t_{j}- t_{j^{\prime}})^{k^{\prime}
+l+1}},\\ j=1,\ldots, n;\ k=1,\ldots,p_j.
\end{array}\end{equation} \begin{dfn}  We shall call
system~(\ref{schl}) the  generalized Schlesinger system.
\end{dfn} Thus we obtain  the first main result of this Thesis:
 \begin{thrm}
\label{schlthrm} System~(\ref{deform}) gives the isoprincipal
 deformation  if,  and only if, $Q_{j,k}(t)$ satisfy
  the  generalized Schlesinger system.
\end{thrm} \begin{rem} \label{firint} \begin{enumerate} \item
Condition (\ref{regt}) is compatible with the generalized
Schlesinger system. Indeed, (\ref{schl}) implies \begin{equation}
\sum_{j=1}^n\d Q_{j,0}=-\sum_{j=1}^n\d Q_{j,0}=0, \end{equation}
\ie $\sum_{j=1}^n Q_{j,0}(t)$ is a first integral of the
generalized Schlesinger system. \item One can check by
straightforward computation that the generalized Schlesinger
system is compatible, and thus completely integrable. Since in the
next section we are going to prove independently a stronger result
(Lemma \ref{solschl}), we do not give the details here.
\end{enumerate} \end{rem}
\section{The Painlev\'e Property}
 Let us consider the Cauchy problem for the
generalized Schlesinger system.  Let
 us fix point $t^0\in{\cal U}$
and, for $j=1,\ldots,n,$ $k=0,\ldots,p_j,$ matrices $Q_{j,k}^0\in
gl (m,{\mathbb C}).$ and consider  system~(\ref{schl}) in ${\cal
U}$ with the initial condition  \begin{equation}\label{inischl}
Q_{j,k}(t^0)=Q_{j,k}^0,\quad j=1,\ldots,n,k=0,\ldots,p_j.
\end{equation}  To begin with, let us assume that $Q_{j,k}^0$
satisfy \begin{equation}\label{reg0} \sum_{j=1}^{n} Q_{j,0}^0=0.
\end{equation} Let ${\cal L}$ be a polydisk, such that
\begin{equation} \ov{\cal U}\subset{\cal L}\subset{\cal S},
\end{equation} and for $j=1,\ldots,n$ let $H_j(x,t^0)\in
GL(m,{\cal O}({\cal L}_j)),$ $M_j(\tilde x)\in GL(m,{\cal
O}(\tilde{\mathbb C}^*))$ be the non-singular and principal
factors of  $Y(\tilde x,t^0)\in GL(m,{\cal O}(\tilde{\cal
R}_{t^0})),$ which satisfies the system
\begin{equation}\label{linsys0}
\der{Y}{x}=\sum_{j=1}^n\sum_{k=0}^{p_{j}} \frac{Q_{j,k}^0}{(x-
t_{j}^0)^{k+1}}Y, \end{equation} with the initial condition
\begin{equation}\label{inisys0} Y(\infty,t^0)\equiv I.
\end{equation} By Lemma~\ref{propfaclmm} $M_j(\tilde x)$ satisfies
in $\tilde{\mathbb C}^*$ linear system (\ref{linm}), with
$J_{j,k}\in gl (m,{\mathbb C})$ are given by
\begin{equation}\label{a0toj} J_{j,k} =\sum_{k^{\prime}=0}^{p_{j}-
k}\oint_{\gamma_j}\frac{ H_{j}^{-1}(x,t^0)
Q_{j,k+k^{\prime}}^0H_{j}(x,t^0)}{(x-
t_j^0)^{k^{\prime}+1}}\frac{\d x}{2\pi i}\end{equation} and
$\gamma_j$ is a homotopically non-trivial simple contour in ${\cal
L}_j\sm\ov{\cal U}_j.$ Let us observe that  the ratio
$M_{j}(\tilde x-t_{j}^0)M_{j}^{-1}(\tilde x- t_{j})$  is a
univalued matrix function in $\Set{{\cal L}_j\sm\ov{\cal
U}_j}\times{\cal U}$ and define \begin{equation}\label{defqxt}
F(x,t)\in GL(m,{\cal O}(\bigcup_{j=1}^n\Set{{\cal L}_j\sm\ov{\cal
U}_j}\times{\cal U}))
\end{equation} by \begin{equation}\label{defqxt2}\begin{array}{rl}
 F(x,t)=H_j(x,t^0)M_{j}(\tilde x-t_{j}^0)&M_{j}^{-1}
(\tilde x-t_{j})H_j^{-1}(x,t^0), \\
 &\quad x\in
{\cal L}_j\sm\ov{\cal U}_j,t\in{\cal U},j=1,\ldots,n.\end{array}
\end{equation} Since \begin{equation}\label{regq} F(x,t^0)\equiv I, \end{equation} it follows from
Lemma~\ref{h3} in Chapter~\ref{intro}
 that there exist a  neighborhood ${\cal V}$ of $t^0$
and matrix functions
\begin{eqnarray} \label{begx+} X^+(x,t)&\in&
GL(m,{\cal O}(\bigcup_j{\cal L}_j\times{\cal V})), \\
\label{begx-} X^-(x,t)&\in &GL(m,{\cal O}(\Set{\ov{\mathbb
C}\sm\bigcup_j\ov{\cal U}_j}\times{\cal V})),
\end{eqnarray} such that \begin{equation}\label{rhbpu} X^-(x,t)F(x,t)=X^+(x,t),\quad
(x,t)\in\bigcup_j\Set{{\cal L}_j\sm\ov{\cal U}_j}\times{\cal V}
\end{equation} and \begin{equation}\label{regx-} X^-(\infty,t)\equiv I.\end{equation} Also, $X^+(x,t)$ and
$(X^+)^{-1}(x,t)$ are meromorphic with respect to  $t$ in ${\cal
U}.$ Now we can formulate the following
\begin{lmm}
\label{solschl} Solution $Q_{j,k}(t)$ of system~(\ref{schl}) with
initial condition (\ref{inischl}), satisfying (\ref{reg0}), is
given in ${\cal U}$ by \begin{equation}\label{xtoa} Q_{j,k}(t)
=\sum_{k^{\prime}=0}^{p_{j}- k}\oint_{\gamma_j}\frac{
X^+(x,t)H_j(x,t^0) J_{j,k+k^{\prime}} H_j^{-
1}(x,t^0)(X^+)^{-1}(x,t)}{(x- t_j)^{k^{\prime}+1}}\frac{\d x}{2\pi
i}. \end{equation}  In particular,  $Q_{j,k}(t)$ are holomorphic
in ${\cal V}$ and meromorphic in ${\cal U}.$ \end{lmm}
\begin{proof} Since  (\ref{regq})
implies \begin{equation}\label{regx+} X^+(x,t^0)\equiv I,
\end{equation}  we can define \begin{equation} H_j(x,t)\in
GL(m,{\cal O}({\cal L}_j\times{\cal V}))
\end{equation} by \begin{equation} H_j(x,t)=X^+(x,t)H_j(x,t^0). \end{equation}
Then (\ref{rhbpu})
implies that there exists \begin{equation} Y(\tilde x,t)\in
GL(m,{\cal O}({\cal R}\cap\psi^{- 1}(\Set{\ov{\mathbb
C}\times{\cal V}}\sm\bigcup_j{\cal T}_j))), \end{equation} which
admits representation (\ref{fact}) in $\tilde{\cal L}_{{\cal
T}_j}\cap\psi^{- 1}(\Set{{\cal L}_j\times{\cal V}}\sm{\cal T}_j)$
for $j=1,\ldots,n$ and the representation \begin{equation}
Y(\tilde x,t)=X^-(x,t)Y(\tilde x,t^0) \end{equation} in
$\psi^{-1}(\Set{\ov{\mathbb C}\sm\cup_j\ov{\cal U}_j}\times{\cal
V}).$ In particular, (\ref{regx-}) means that $Y(\infty, t)\equiv
I$  for $t$ in ${\cal V}.$ According to Lemma~\ref{propfaclmm},
$Y(\tilde x,t)$ satisfies system~(\ref{deform}), where
$Q_{j,k}(t)\in gl(m,{\cal O}({\cal V}))$ are given by (\ref{xtoa})
and satisfy (\ref{regt}), (\ref{inischl}). By our construction,
this system gives the isoprincipal deformation in ${\cal V}.$
Hence, by Theorem~\ref{schl}, $Q_{j,k}(t)$ satisfy
system~(\ref{schl}) in ${\cal V}.$ Finally, since $X^+(x,t)$ and
$(X^+)^{-1}(x,t)$ are meromorphic with respect to $t$ in ${\cal
U},$ we conclude that $Q_{j,k}(t),$ defined by (\ref{xtoa}), give
in ${\cal U}$ meromorphic solution to system~(\ref{schl}) with
initial condition (\ref{inischl}).
\end{proof}
 We would
like to note here that there is no loss of generality in
assumption (\ref{reg0}).  Indeed, we can always fix  $a
\in{\mathbb C}\sm\cup_j\ov{\cal L}_j$ and define map $\mu:{\cal
U}\mapsto{\mathbb C}^{n+1}$ by \begin{equation}
 \mu(t_1,\ldots,t_n)=(\frac{1}{t_1-
a},\ldots,\frac{1}{t_n-a},0).
\end{equation}
 Furthermore, for
$j=1,\ldots,n+1,$ $k=0,\ldots,p_j,$ $p_{n+1}=0,$ we can consider
solution $Q_{j,k}^{\prime}(\tau)$ of the generalized Schlesinger
system in $n+1$ variables $\tau=(\tau_1,\ldots, \tau_{n+1})$ with
the initial condition \begin{equation}
 \begin{array}{l}
Q_{j,k}^{\prime}(\mu(t^0))=Q_{j,k}^0,\quad j=1,\ldots,n, \\
Q_{n+1,0}^{\prime}(\mu(t^0))=-\sum_{j=1}^{n}Q_{j,0}^0, \end{array}
\end{equation}
 so that   \begin{equation}
 \sum_{j=1}^{n+1}Q_{j,0}^{\prime}(\mu(t^0))=0.
\end{equation} According to Lemma~\ref{solschl}, $Q_{n+1,0}^{\prime}(\tau)$
is of the form \begin{equation}\label{an0}
Q_{n+1,0}^{\prime}(\tau)=X(\tau)Q_{n+1,0}^{\prime}(\mu(t^0))X^{-1}(\tau),
\end{equation} where meromorphic matrix function $X(\tau)$  satisfies
\begin{equation}\label{regxt0} X(\mu(t^0))=I \end{equation} and, in view of
Lemma~\ref{schlmm}, \begin{equation}\label{sysx}
\dder{X}{\tau_j}X^{-1}(\tau)= -
\sum_{k=0}^{p_j}\frac{Q_{j,k}^{\prime}(\tau)}{(\tau_{n+1}-
\tau_j)^{k+1}},\quad j=1,\ldots,n. \end{equation} Now it can be
checked by straightforward computation that $Q_{j,k}(t),$  defined
by \begin{equation} Q_{j,k}(t)=\left (X^{-1}Q_{j,k}^{\prime}
X\right )\circ\mu(t),\quad j=1,\ldots,n; k=0,\ldots,p_j,
\end{equation} satisfy system~\ref{schl} with initial condition
(\ref{inischl}).

Thus solution of system~(\ref{schl}) with arbitrary initial values
at $t^0$ is meromorphic in ${\cal U}.$ Since our choice of ${\cal
U}$ in the first place was limited only by the assumption
$\ov{\cal U}\subset{\cal S},$ this solution can be continued into
universal covering space $\tilde {\cal S}$ over ${\cal S},$ and
we obtain the second main result of this Thesis:
\begin{thrm} \label{tm}  Any solution of the
generalized Schlesinger system~(\ref{schl}) is meromorphic in
$\tilde{\cal S}.$ \end{thrm} Theorem~\ref{tm} implies, in
particular, that the generalized Schlesinger system enjoys the
{\em Painlev\'e property,} that is, it has no movable critical
points.
\chapter{Fuchsian Systems} \label{Fuchs}
 \section{Non-Resonant Fuchsian Systems}
 \begin{dfn} A linear differential system of
the form \begin{equation}\label{Ful}
\der{Y}{x}=\sum_{j=1}^n\frac{Q_{j}}{x-t_{j}}Y,
\end{equation} where $Q_j\in gl(m,{\mathbb C}),$ is called a Fuchsian system.
\end{dfn}  System~(\ref{Ful}) is a special case of
system~(\ref{linsys}) with \begin{equation}\label{Fu}
p_1=\ldots=p_n=0,\quad Q_{j,0}=Q_j. \end{equation} As before, we
assume  that $x=\infty$ is a regular point, \ie
\begin{equation}\label{Fureg} \sum_{j=1}^nQ_j=0,
\end{equation} and set the initial condition \begin{equation}\label{inif}
Y(\infty)=I. \end{equation} Then we can apply the theory developed
in the previous Chapter to the Fuchsian systems.

To begin with, we consider the following \begin{dfn} Fuchsian
system~(\ref{Ful}) is called non-resonant if for $j=1,\ldots,n$ no
two eigenvalues of $Q_j\in gl(m,{\mathbb C})$ differ by an
integer.
\end{dfn}
 It is well-known\footnote{See any
textbook on Fuchsian systems, such as \cite{codlev}.} that
for the non-resonant Fuchsian systems  most of the
information about the monodromy representation and the local
factorization of the fundamental solution can be obtained very
easily:
 \begin{lmm} \label{Fufaclmm} Let system~(\ref{Ful}) be
non-resonant. Then for $j=1,\ldots,n$ $Y$ has at $t_j$ the
principal factor of the form \begin{equation}\label{nonprin}
M_j(\tilde x)=\tilde x^{J_j}K_j, \end{equation} where $J_j$ is the
Jordan form of $Q_j,$ and $K_j\in GL(m,{\mathbb C}).$ In
particular, the corresponding monodromy generator $\Phi(\sigma_j)$
is given by \begin{equation}\label{nonmon}
\Phi(\sigma_j)=K_j^{-1}e^{-2\pi iJ_j}K_j.
\end{equation}
\end{lmm}
\begin{proof} We give a sketch of the proof, adapted
from \cite{codlev}.
 Let us fix $j.$ Considering, if necessary, the
transformation \begin{equation} Y^{\prime}=TY, \end{equation}
where \begin{equation} TQ_jT^{-1}=J_j, \end{equation}  we can
assume that \begin{equation} Q_j=J_j.
\end{equation} In view of Lemma~\ref{propfaclmm}, it suffices to show that
the system \begin{equation}\label{holsysFu}
\der{H_j}{x}=\frac{[J_j,H_j(x)]}{x-
t_j}+\sum_{j^{\prime}\not=j}\frac{Q_{j^{\prime}}}
{x-t_{j^{\prime}}}H_j(x) \end{equation} has solution $H_j(x),$
holomorphic in a neighborhood of $t_j$ and normalized by
\begin{equation} H_j(t_j)=I. \end{equation} Considering the Taylor
expansions
\begin{eqnarray}
H_j(x)&=&I+\sum_{k=1}^{\infty}H_{j,k}(x-t_j)^k, \\
\sum_{j^{\prime}\not=j}\frac{Q_{j^{\prime}}} {x-
t_{j^{\prime}}}&=&\sum_{k=0}^{\infty}R_{j,k}(x-t_j)^k,
\end{eqnarray}
we obtain the recurrence relation for $H_{j,k}:$
\begin{equation}\label{recur}
(k+1)H_{j,k+1}-[J_j,H_{j,k+1}]=\sum_{k^{\prime}=0}^kR_{j,k^{\prime}}H_{j,k-
k^{\prime}}. \end{equation} The left-hand side of (\ref{recur})
can be considered as an expression for a linear operator, acting
on $gl(m,{\mathbb C}).$ The non-resonance condition guarantees
that this operator is invertible, \ie the recurrence relation
(\ref{recur}) is solvable. The convergence of the obtained Taylor
series for $H_j(x)$ can be
 established by the Cauchy majorization method. This completes the
proof.
\end{proof}
Substituting (\ref{Fu}) into (\ref{schl}), we observe that in the
case of Fuchsian systems the  generalized Schlesinger system takes
the form \begin{equation}\label{Sch} \d Q_{j}
=\sum_{j^{\prime}\not=j} \left ]Q_{j^{\prime}}, Q_{j}\right
]\frac{\d t_{j}-\d t_j^{\prime}} {t_j-t_j^{\prime}},\quad
j=1,\ldots,n. \end{equation} System~(\ref{Sch}) was introduced by
L.~Schlesinger in \cite{schl} and is known as {\em the Schlesinger
system.} This explains our term "generalized" for
system~(\ref{schl}).
 In order to consider the
implications of Lemma~\ref{Fufaclmm} for the Schlesinger system,
let us recall  the notation \begin{equation}\label{defset} {\cal
S}={\mathbb
C}^n\sm\bigcup_{j^{\prime}\not=j^{\prime\prime}}\Set{t=(t_1,\ldots,t_n):t_{j^{\prime}
} =t_{j^{\prime\prime}}}, \end{equation} fix point $t^0\in {\cal
S},$ matrices $Q_j^0\in gl(m,{\mathbb C}),$ such that
\begin{equation}\label{Fureg0}\sum_{j=1}^n Q_j^0=0, \end{equation} and assume that
the Fuchsian system \begin{equation}\label{Ful0}
\dder{Y}{x}=\sum_{j=1}^n\frac{Q_{j}^0}{x-t_{j}^0}Y \end{equation}
is non-resonant. Now we  consider  an isomonodromic deformation of
system~(\ref{Ful0}). Let polydisk ${\cal U}$  and  matrix
functions $Q_1(t),\ldots,Q_n(t)\in gl(m,{\cal O}({\cal U}))$ be
such that
\begin{eqnarray}
&& t^0\in{\cal U}\subset\ov{\cal U}\subset{\cal S},\\
\label{Furegt} && \sum_{j=1}^n
Q_j(t)\equiv 0, \\ \label{iniSch} && Q_j(t^0)=Q_j^0,
\end{eqnarray}
 and system
\begin{equation}\label{Fudeform} \dder{Y}{x}=\sum_{j=1}^n\frac{Q_j(t)}{x-t_{j}}Y,
\end{equation} with initial condition \begin{equation}\label{Fuinidef} Y(\infty,t)
\equiv  I \end{equation} gives an isomonodromic deformation in
${\cal U}.$ Then, by continuity, system~(\ref{Fudeform}) is
non-resonant for each fixed $t$ in a neighborhood of $t^0.$ Since
for $j=1,\ldots,n$ the monodromy generator $\Phi(\sigma_j)$ is
independent of $t,$ Lemma~\ref{Fufaclmm} implies that $Q_j(t)$ has
 Jordan form $J_j,$ independent of $t$ in the whole of ${\cal
U}.$ Since matrix $K_j$ in expression~(\ref{nonmon}) is determined
by matrices $\Phi(\sigma_j)$ and $J_j$ up to multiplication from
the left by a non-degenerate  matrix, commuting with $J_j,$  we
can observe in view of Remark~\ref{facrem} that the principal
factors of $Y(\tilde x,t)$ can be chosen to be independent of $t.$
Hence, system~(\ref{Fudeform}) actually gives the  isoprincipal
deformation in ${\cal U}.$  On the other hand, we know (see
Remark~\ref{schliso}) that the isoprincipal deformation is
isomonodromic. Thus Theorem~\ref{schlthrm} implies the following
result, originally proved by  L.~Schlesinger in \cite{schl}:
\begin{thrm}[Schlesinger] \label{Schthrm}
Assume that Fuchsian
system~(\ref{Ful0}) is non-resonant. Then $Q_j(t)$ satisfy
system~(\ref{Sch}) with initial condition (\ref{iniSch}) if, and
only if, system~(\ref{Fudeform}) gives the isomonodromic
deformation of Fuchsian system~(\ref{Ful0}).
\end{thrm}
We would like to note that without  the assumption of
non-resonance the conclusion of Theorem~\ref{Schthrm} fails in one
direction. There exist isomonodromic deformations of resonant
Fuchsian systems, which are not governed by the Schlesinger
system. In the subsequent sections of this Chapter we shall
present V.~E.~Katsnelson's example,
 concerning
rational matrix functions in general position, illustrating this
phenomenon.

Theorem~\ref{schlthrm}, on the other hand, is applicable in the resonant
 case as well,\ie in order  to solve the Cauchy problem for the Schlesinger
system one has to construct the isoprincipal deformation. We shall
see that in the case of rational matrix functions in general
position this approach allows to solve the Cauchy problem for the
Schlesinger system explicitly. In this context we would like to
mention a generalization of Lemma~\ref{Fufaclmm} due to
A.~H.~M.~Levelt (see \cite{levelt}). It states that in general
fundamental solution $Y(x)$ of Fuchsian system~(\ref{Ful}) has at
$t_j$ the principal factor of the form \begin{equation} M_j(\tilde
x)=x^{D_j}\tilde x^{E_j}K_j, \end{equation} where $D_j$ is a
diagonal matrix with integer entries, $E_j$ is an upper-triangular
matrix, whose spectrum coincides modulo ${\mathbb Z}$ with the
spectrum of $Q_j$ and consists of eigenvalues, pairwise distinct
modulo ${\mathbb Z},$
 and $K_j\in GL(m,{\mathbb C}).$

We would also like to note that Theorem~\ref{tm} implies
that any solution of the
 Schlesinger system~(\ref{Sch}) is meromorphic in the universal
covering over ${\cal S}.$ This result was originally proved by
T.~Miwa in \cite{miwa} under the assumption of non-resonance,
which now can be removed.
 \section{Rational Matrix Functions in General Position}
Now our goal is to deal with the class of Fuchsian systems, whose
fundamental solutions are rational matrix functions in
general position. General theory of such systems
was discussed in depth in V.~E.~Katsnelson's works \cite{katze1},
\cite{katze}. We follow these sources in the presentation below.
 \begin{dfn} Rational square matrix
function $Y(x)$ is said to be  a rational matrix function in
general position, if: \begin{enumerate} \item $\det{Y(x)}$ does
not vanish identically; \item the polar set of $Y(x)$  and the
polar set  of $Y^{-1}(x)$  do not intersect; \item all  the poles
of $Y(x)$ and $Y^{-1}(x)$ are simple; \item all  the residues of
$Y(x)$ and $Y^{-1}(x)$ are  matrices of rank one; \item both
$Y(x)$ and $Y^{-1}(x)$ are holomorphic at $\infty.$
\end{enumerate} \end{dfn} Let $Y(x)$ be an $m\times m$ rational
matrix function in general position. Without loss of generality,
we assume \begin{equation}\label{inirat} Y(\infty)=I.
\end{equation}  Let us order somehow the poles of $Y(x)$
(respectively, the poles of $Y^{-1}(x)$) and denote them by
$t_1,\ldots,t_s$ (respectively, by
$t_{s+1},\ldots,t_{s+s^{\prime}}$). Since, by definition, all  the
poles of $Y(x)$ and $Y^{-1}(x)$ are simple and all the residues of
$Y(x)$ and $Y^{-1}(x)$ have rank one, it follows from
(\ref{inirat}) that \begin{equation}
\det{Y(x)}=\frac{\prod_{j=1}^{s^{\prime}}(x- t_{s+j})}{
\prod_{j=1}^s(x-t_j)}, \end{equation} and, in particular, that
$s=s^{\prime}.$ Thus the polar sets of $Y(x)$ and $Y^{-1}(x)$ have
the same cardinality.  We shall call the ordered sets
\begin{eqnarray}
\label{poleset} {\cal P}&=&\Set{t_1,\ldots,t_s}, \\
\label{zeroset} {\cal Z}&=&\Set{t_{s+1},\ldots,t_{2s}},
\end{eqnarray}
 respectively,
the {\em pole} and {\em zero sets} of $Y(x).$ We shall also
associate with these sets the diagonal matrices
\begin{eqnarray}
A_{{\cal P}}&=&\diag{t_1,\ldots,t_s}, \label{ap} \\
A_{{\cal Z}}&=&\diag{t_{s+1},\ldots,t_{2s}}, \label{an}
\end{eqnarray}
which we shall call, respectively, the {\em pole} and {\em zero
matrices} of $Y(x).$ Since for $j=1,\ldots,s,$ \begin{equation}
\rank{\res{Y}{x= t_{j}}}=1, \end{equation} there exist $m\times 1$
matrix $c_{j}$ and $1\times m$ matrix $b_{j},$ such that
\begin{equation}\label{polesem} \res{Y}{x= t_{j}}=c_{j}b_{j}. \end{equation} We shall
call matrices $c_1,\ldots,c_s$ (respectively, $b_1,\ldots,b_s$)
the {\em left} (respectively, {\em right}) {\em  pole
semi-residues} of $Y(x).$ Analogously, we shall call $m\times 1$
matrices $c_{s+1}, \ldots,c_{2s}$ and $1\times m$ matrices
$b_{s+1}, \ldots,b_{2s},$ such that
  for $j=s+1,\ldots,2s,$
\begin{equation}\label{zerosem} \res{Y^{- 1}}{x=t_{j}}=c_{j}b_{j}, \end{equation} the
{\em left} and {\em right zero semi-residues} of $Y(x).$ From the
left semi-residues of $Y(x)$ we construct  two $m\times s$
matrices \begin{equation}\label{cp} C_{{\cal P}}=\left (c_1\cdots
c_s\right ), \ C_{{\cal Z}}=\left (c_{s+1}\cdots c_{2s}\right ),
\end{equation} which are called, respectively, the {\em left pole}
and {\em zero semi-residual matrices} of $Y(x).$ From the right
semi-residues of $Y(x)$ we construct  two $s\times m$ matrices
\begin{equation}\label{bn}
B_{{\cal P}}=\left (\begin{array}{c}b_1 \\ \vdots \\
b_s\end{array}\right ), \
 B_{{\cal Z}}=\left (\begin{array}{c}b_{s+1} \\ \vdots \\
b_{2s}\end{array}\right ), \end{equation} which are  called,
respectively, the
 {\em right pole} and {\em zero  semi-residual
matrices} of $Y(x).$
\begin{rem}
\label{amsemres}  Note that semi-residual matrices $C_{{\cal
P}},B_{{\cal P}},C_{{\cal Z}},B_{{\cal Z}}$ of $Y(x)$ are defined
up to the transformation \begin{equation} \left\{\begin{array}{l}
C_{{\cal P}}^{\prime}=C_{{\cal P}}E_{{\cal P}},\\
B_{{\cal P}}^{\prime}=E_{{\cal P}}^{-1}B_{{\cal P}},\\ C_{{\cal Z}}^{\prime}=C_{{\cal Z}}E_{{\cal Z}},\\
B_{{\cal Z}}^{\prime}=E_{{\cal Z}}^{-1}B_{{\cal Z}}
\end{array}\right.\end{equation} for arbitrary diagonal $E_{{\cal
P}},E_{{\cal Z}}\in GL(s,{\mathbb C}).$
\end{rem}
In terms of the semi-residual matrices  we
can write down the following   {\em realizations} of $Y(x)$ and $Y^{-1}(x):$
\begin{eqnarray} Y(x)&=&I+C_{{\cal P}}\left (xI-A_{{\cal P}}\right )^{-1} B_{{\cal P}},
\label{repy1} \\ Y^{-1}(x) &=&I+ C_{{\cal Z}}\left (xI-A_{{\cal
Z}}\right )^{-1} B_{{\cal Z}}. \label{repy2} \end{eqnarray} In
order to explain, what kind of relations the identity
\begin{equation} Y(x)Y^{-1}(x)= Y^{-1}(x)Y(x)= I \end{equation}
imposes on the semi-residual matrices of $Y(x),$
  we  formulate the following version
of a result, proved by I.~Gohberg, M.~A.~Kaashoek, L.~Lerer
and L.~Rodman in \cite{gohberg}:
 \begin{thrm}[Gohberg, Kaashoek, Lerer, Rodman]
 \label{t1.1}
\begin{enumerate}
\item Let $Y(x)$ be a rational matrix function in general
position, normalized at $\infty$ by (\ref{inirat}), with  pole and
zero matrices $A_{{\cal P}},A_{{\cal Z}}$ and left pole, left
zero, right pole, right zero
 semi-residual matrices  $C_{{\cal P}},C_{{\cal Z}},
B_{{\cal P}}, B_{{\cal Z}},$  respectively. Let $S_{{\cal P}{\cal
Z}}$ be the solution of the matricial equation
\begin{equation}\label{lypn} A_{{\cal P}}S_{{\cal P}{\cal Z}}-
S_{{\cal P}{\cal Z}}A_{{\cal Z}}=B_{{\cal P}}C_{{\cal Z}}
 \end{equation}
and let $S_{{\cal Z}{\cal P}}$ be the   solution
 of the matricial equation
\begin{equation}\label{lynp} A_{{\cal Z}}S_{{\cal Z}{\cal P}}- S_{{\cal Z}{\cal
P}}A_{{\cal P}}=B_{{\cal Z}}C_{{\cal P}}.
 \end{equation}
Then matrices $S_{{\cal P}{\cal Z}}$ and $S_{{\cal Z}{\cal P}}$
are mutually inverse: \begin{equation}\label{incore} S_{{\cal
P}{\cal Z}}=S_{{\cal Z}{\cal P}}^{-1}, \end{equation} and the
following relations hold true:
\begin{eqnarray}
\label{exe1}
C_{{\cal P}} S_{{\cal P}{\cal Z}}&=& C_{{\cal Z}},\\
\label{exe2} S_{{\cal P}{\cal Z}}B_{{\cal Z}}&=& -B_{{\cal P}}.
 \end{eqnarray}
\item
 Let
${\cal P},{\cal Z}$  be two  disjoint finite  subsets of ${\mathbb
C}$ of the same cardinality $s$ and let $A_{{\cal P}},A_{{\cal
Z}}$ be the associated diagonal matrices. Let $C_{{\cal
Z}},B_{{\cal P}}$ (respectively, $C_{{\cal P}},B_{{\cal Z}}$) be
an $m\times s$ matrix and an $s\times m$ matrix. Assume that
$s\times s$ matrix solution $S_{{\cal P}{\cal Z}}$ (respectively,
$S_{{\cal Z}{\cal P}}$) of equation (\ref{lypn}) (respectively,
(\ref{lynp})) is invertible.
 Then there exists  unique $m\times m$
rational matrix function in general position $Y(x),$ normalized at
$\infty$ by (\ref{inirat}), with  pole and zero matrices $A_{{\cal
P}},A_{{\cal Z}}$ and left zero, right pole (respectively, left
pole, right zero)
 semi-residual matrices  $C_{{\cal Z}},
B_{{\cal P}}$  (respectively, $C_{{\cal P}},B_{{\cal Z}}$).
\end{enumerate}
 \end{thrm}
\begin{proof} \begin{enumerate}
\item   Substituting expressions (\ref{repy1}), (\ref{repy2}) and
(\ref{lypn}), rewritten in the form \begin{equation} B_{{\cal
P}}C_{{\cal Z}}=S_{{\cal P}{\cal Z}}(xI-A_{{\cal Z}})-(xI-A_{{\cal
P}})S_{{\cal P}{\cal Z}}, \end{equation} into the identity
$Y(x)Y^{-1}(x)= I ,$ we obtain \begin{equation} C_{{\cal P}}\left
(xI-A_{{\cal P}}\right )^{-1}\left (B_{{\cal P}}+S_{{\cal P}{\cal
Z}}B_{{\cal Z}}\right )=
 =\left (C_{{\cal Z}}-C_{{\cal P}}S_{{\cal P}{\cal Z}}\right )\left (xI-A_{{\cal P}}\right )^{-1}B_{{\cal P}}. \end{equation}
Since the spectra of $A_{{\cal P}}$ and $A_{{\cal Z}}$ are
disjoint, this means that relations (\ref{exe1}) and (\ref{exe2})
hold true.
 Analogously,  we can derive from equation (\ref{lynp}) and the
identity $Y^{-1}(x)Y(x)= I$ the following relations:
\begin{eqnarray}
\label{exe11}
C_{{\cal Z}} S_{{\cal Z}{\cal P}}&=& C_{{\cal P}},\\
\label{exe21} S_{{\cal Z}{\cal P}}B_{{\cal P}}&=& -B_{{\cal Z}}.
 \end{eqnarray}
It follows from (\ref{lypn}) and (\ref{exe11}) that
\begin{equation} A_{{\cal P}}S_{{\cal P}{\cal Z}}S_{{\cal Z}{\cal
P}}-S_{{\cal P}{\cal Z}}A_{{\cal Z}}S_{{\cal Z}{\cal P}}=B_{{\cal
P}}C_{{\cal P}}.
\end{equation}
 On the other hand, it follows
from (\ref{lynp}) and (\ref{exe2}) that
 \begin{equation} S_{{\cal P}{\cal Z}}A_{{\cal Z}}S_{{\cal Z}{\cal P}}-
 S_{{\cal P}{\cal Z}}S_{{\cal Z}{\cal P}}A_{{\cal P}}=-B_{{\cal P}}C_{{\cal P}}.
\end{equation} Thus \begin{equation} \left  [A_{{\cal P}}, S_{{\cal P}{\cal Z}}S_{{\cal
Z}{\cal P}}\right  ]=0, \end{equation} and, since all eigenvalues
of $A_{{\cal P}}$ are of multiplicity $1,$ matrix $S_{{\cal
P}{\cal Z}}S_{{\cal Z}{\cal P}}$ is diagonal. Since (\ref{exe1})
and (\ref{exe11}) imply \begin{equation} C_{{\cal P}}S_{{\cal
P}{\cal Z}}S_{{\cal Z}{\cal P}}=C_{{\cal P}}, \end{equation} and
since, by definition, left pole semi-residual matrix $C_{{\cal
P}}$ has no zero columns, we conclude that \begin{equation}
S_{{\cal P}{\cal Z}}S_{{\cal Z}{\cal P}}=I. \end{equation} \item
The proof is absolutely analogous for both versions, so we assume
that matrices $C_{{\cal Z}}, B_{{\cal P}}$ are given and that
matrix $S_{{\cal P}{\cal Z}},$ satisfying equation (\ref{lypn}),
is invertible. If rational matrix function in general position
$Y(x),$ normalized at $\infty$ by (\ref{inirat}), with  pole and
zero matrices $A_{{\cal P}},A_{{\cal Z}}$ and left zero, right
pole
 semi-residual matrices  $C_{{\cal Z}},
B_{{\cal P}}$  exists, then, in view of part 1) of the Theorem, it
admits realizations (\ref{repy1}), (\ref{repy2}) with left pole,
right zero  semi-residual matrices
 $C_{{\cal P}},B_{{\cal Z}}$  given by
\begin{eqnarray}
\label{exe12}
C_{{\cal P}} &=& C_{{\cal Z}}S_{{\cal P}{\cal Z}}^{-1},\\
\label{exe22} B_{{\cal Z}}&=& -S_{{\cal P}{\cal Z}}^{-1}B_{{\cal
P}}.
 \end{eqnarray}
Hence it is  uniquely determined. In order to prove the existence
of such matrix function $Y(x),$ it suffices
to verify that expressions (\ref{repy1}), (\ref{repy2}) agree with
each other, \ie that the
identity  $Y(x)Y^{-1}(x)= I $ holds true.
This can be done by straightforward computation, utilizing (\ref{lypn})
in the same way as above.
\end{enumerate}
 \end{proof}
\begin{rem}
\label{gohberg}
\begin{enumerate}
\item In principle, the existence and uniqueness of solution for
equations~(\ref{lypn}), (\ref{lynp}) follows from the
disjointedness of the spectra of matrices $A_{{\cal P}}$ and
$A_{{\cal Z}}.$ However, since these matrices are diagonal, we can
give the explicit solutions of these equations, with notations
(\ref{ap}), (\ref{an}), (\ref{cp}), (\ref{bn}):
\begin{eqnarray}
\label{spn} S_{{\cal P}{\cal Z}}&=&\left
(\frac{b_{\alpha}c_{s+\beta}}{t_{\alpha}-t_{s+\beta}}\right )_{
\alpha,\beta=1}^s,\\
\label{snp} S_{{\cal Z}{\cal P}}&=&\left
(\frac{b_{s+\alpha}c_{\beta}}{t_{s+\alpha}-t_{\beta}}\right )_{
\alpha,\beta=1}^s.
\end{eqnarray}
If $Y(x)$ is a rational matrix function in general position, we
shall call matrices $S_{{\cal P}{\cal Z}}$ and $S_{{\cal Z}{\cal
P}},$ constructed in (\ref{spn}), (\ref{snp}) from the
semi-residues, the poles and the zeroes of $Y(x),$ the {\em
pole-zero} and {\em zero-pole core matrices} of $Y(x).$ \item
Actually, the result proved in \cite{gohberg} is more general than
Theorem~\ref{t1.1}. It concerns invertible rational matrix
functions, not necessarily in general position. Such matrix
functions also have realizations of form (\ref{repy1}),
(\ref{repy2}). However, in general, $s\times s$ matrices $A_{{\cal
P}}, A_{{\cal Z}}$ in these
 expressions
need not be diagonal and their spectra need not be disjoint. If
realization (\ref{repy1}) is minimal, \ie dimension $s$ is minimal
possible, then the pair $(A_{{\cal P}},B_{{\cal P}})$ is called a
pole pair of $Y(x).$ If realization (\ref{repy2}) is minimal, then
the pair $(C_{{\cal Z}},A_{{\cal Z}})$ is called a zero pair of
$Y(x).$
 Under the assumption that the pairs $(C_{{\cal Z}},A_{{\cal Z}})$ and $(A_{{\cal P}},B_{{\cal P}})$
are minimal, \i.e. that
\begin{eqnarray}
\bigcap_{j=0}^{s-1}\Ker{C_{{\cal Z}}A_{{\cal Z}}^j}&=&0,\\
\sum_{j=0}^{s-1}\Im{A_{{\cal P}}^j}B_{{\cal P}}&=&{\mathbb C}^s,
\end{eqnarray}
Gohberg-Kaashoek-Lerer-Rodman theorem asserts that $(C_{{\cal
Z}},A_{{\cal Z}})$ and $(A_{{\cal P}},B_{{\cal P}})$ are a zero
pair and a pole pair of some rational matrix function $Y(x)$
 if, and only
if, equation (\ref{lypn}) has an invertible solution.
There is one-to-one correspondence between such rational matrix
functions and invertible solutions of equation (\ref{lypn}),
given by (\ref{exe1}), (\ref{exe2}).
\end{enumerate}
\end{rem}
From Theorem~\ref{t1.1} we can derive several results,
illustrating the connection between rational matrix functions in
general position and Fuchsian systems. For brevity we shall use
the following notation: $I_j$ is the diagonal matrix
\begin{equation} I_j=\diag{0,\ldots,0,1,0,\ldots,0},
\end{equation} whose only non-zero element is in the $j$-th
position on the diagonal.
\begin{lmm} \label{Fuchsrat}
Let $Y(x)$ be a rational matrix function in general position,
normalized at $\infty$ by (\ref{inirat}), with  pole and zero
matrices $A_{{\cal P}},A_{{\cal Z}}$ of form  (\ref{ap}),
(\ref{an}), and  left pole,  right zero
 semi-residual matrices  $C_{{\cal P}},
B_{{\cal Z}}.$   Let $S_{{\cal Z}{\cal P}}$ be the zero-pole core
matrix of $Y(x).$ Then $Y(x)$ is a fundamental solution of  the
Fuchsian system \begin{equation}\label{Furat}
\der{Y}{x}=\sum_{j=1}^{2s}\frac{Q_{j}}{x-t_{j}}Y, \end{equation}
where for $j=1,\ldots,2s$ $Q_j$ is given by
\begin{equation}\label{rataj} Q_{j}=\left\{
\begin{array}{l@{\quad}l} C_{{\cal P}}I_jS_{{\cal Z}{\cal P}}^{-1}\left
(t_{j}I-A_{{\cal Z}}\right )^{-1} B_{{\cal Z}},&j=1,\ldots,s,\\
C_{{\cal P}}\left (t_{j}I-A_{{\cal P}}\right )^{- 1}S_{{\cal
Z}{\cal P}}^{-1}I_{j-s} B_{{\cal Z}},&j=s+1,\ldots,2s. \end{array}
\right. \end{equation}
\end{lmm}
\begin{proof}
In the manner completely analogous
to the proof of Theorem~\ref{t1.1}, we can derive from
 (\ref{repy1}), (\ref{repy2})  and (\ref{lypn}) the following expression for the
logarithmic derivative of $Y(x):$ \begin{equation}\label{qx}
\begin{array}{rl} Q(x)&=\der{Y}{x}Y^{-1}(x) \\ &=C_{{\cal P}}\left
(xI-A_{\cal P}\right )^{-1} S_{{\cal P}{\cal Z}}\left (xI-A_{{\cal
Z}}\right )^{-1} B_{{\cal Z}}.
\end{array} \end{equation}
 In order
to obtain (\ref{rataj}), it remains only to substitute (\ref{incore})
into (\ref{qx}) and to compute the residue
of $Q(x)$ at $t_j.$
\end{proof}
In view of Lemma~\ref{Fuchsrat}, it makes sense to consider
the principal factors of a rational matrix function in general
position.
\begin{lmm}
\label{prinrat} Let ${\cal P},{\cal Z}$  be two  disjoint finite
subsets of ${\mathbb C}$ of the same cardinality $s,$ with
notation (\ref{poleset}), \ref{zeroset}), and let
 $Q_j\in gl(m,{\mathbb C})$ for $j=1,\ldots,2s.$ Assume that
\begin{equation} \sum_{j=1}^{2s}Q_j=0 \end{equation} and let $Y(x)$ be the
fundamental solution of Fuchsian system (\ref{Furat}) with initial
condition (\ref{inirat}). Then $Y(x)$ is a rational matrix
function in general position with pole and zero sets ${\cal
P},{\cal Z}$ if, and only if, $Y(x)$ has the principal factors of
the form \begin{equation}\label{defprinrat} M_{j}(x)= \left\{
\begin{array}{l@{\quad}l} x^{- I_1}K_{j},&j=1,\ldots,s,\\
x^{I_1}K_{j},&j=s+1,\ldots,2s,
\end{array} \right. \end{equation} where for $j=1,\ldots,2s$ $K_{j}\in Gl(m,{\mathbb C}).$
In this case for $j=1,\ldots,s$ the first row
of matrix $K_j$ is the right pole semi-residue of $Y(x)$
at $t_j$ and  for $j=s+1,\ldots,2s$ the first column
of matrix $K_j^{-1}$ is the left zero semi-residue of $Y(x)$
at $t_j.$
\end{lmm}
\begin{proof}
Let us assume first that $Y(x)$ is a rational matrix function in
general position with pole and zero sets ${\cal P},{\cal Z}.$ Let
$b_1,\ldots,b_2s,c_1,\ldots,c_{2s}$ be the right pole, right zero,
left pole,  left zero semi-residues of $Y(x)$ and
 let us denote by
$i_j$  the $j$-th vector-column of the standard basis. Now for
$j=1,\ldots,2s$ let $K_{j}\in Gl(m,{\mathbb C})$ satisfy
\begin{eqnarray}
\label{bj}
i_1^*K_j=b_j,&\quad&j=1,\ldots,s,\\
\label{cj}
K_j^{-1}i_1=c_j,&\quad&j=s+1,\ldots,2s.
\end{eqnarray}
We are going to prove that for $j=1,\ldots,2s$ matrix function
$M_j(x),$ defined by (\ref{defprinrat}) is the principal factor of
$Y(x)$ at $t_j.$ By definition of the principal factors, we have
to show for $j=1,\ldots,2s$ that matrix function \begin{equation}
H_j(x)=Y(x)M_j^{-1}(x-t_j). \end{equation}  is holomorphic and
invertible in a neighborhood of $t_j.$ Both  $H_j$ and $H_j^{-1}$
have at most a simple pole at $t_j,$ so we only have to compute
the appropriate residues. For $j=1,\ldots,s$ we have
\begin{equation}
\res{H_j}{x=t_j}=c_jb_jK_j^{-1}(I-I_1)=c_ji_1^*(I-I_1)=0,
\end{equation} and also, since \begin{equation}
\res{YY^{-1}}{x=t_j}=c_jb_jY^{- 1}(t_j)=0, \end{equation} we
obtain \begin{equation} \res{H_j^{-1}}{x=t_j}=I_1K_jY^{-
1}(t_j)=i_1b_jY^{-1}(t_j)=0. \end{equation} For $j=s+1,\ldots,2s$
the computations are absolutely analogous.

Conversely, let us assume that for $j=1,\ldots,2s$ $Y(x)$ admits
in a neighborhood of $t_j$ the factorization \begin{equation}
Y(x)=H_j(x)M_j(x-t_j), \end{equation} where matrix function
$H_j(x)$ is holomorphic and invertible in a neighborhood of $t_j$
and matrix function $M_j(x)$ is given by (\ref{defprinrat}). Then
$Y(x)$ is an invertible rational matrix function, ${\cal P}$ is
the polar set of $Y(x)$ and ${\cal Z}$ is the polar set of
$Y^{-1}(x).$ Moreover,
\begin{eqnarray}
\res{Y}{x= t_{j}}=c_{j}b_{j},&\quad&j=1,\ldots,s,\\
\res{Y^{-1}}{x= t_{j}}=c_{j}h_{j},&\quad&j=s+1,\ldots,2s,
\end{eqnarray}
where $c_1,\ldots, c_{2s}$ and $b_1,\ldots,b_{2s}$ are, respectively,
$m\times 1$ and $1\times m$ matrices, given by
\begin{eqnarray}
c_j&=&\left\{\begin{array}{l@{\quad}l}
H_j(t_j)i_1,&j=1,\ldots,s,\\
K_j^{-1}i_1,&j=s+1,\ldots,2s,
\end{array}\right.\\
b_j&=&\left\{\begin{array}{l@{\quad}l}
i_1^*K_j,&j=1,\ldots,s,\\
i_1^*H_j^{-1}(t_j),&j=s+1,\ldots,2s.
\end{array}\right.\end{eqnarray}
Hence, $Y(x)$ is a rational matrix function in general
position and $b_1,\ldots,b_s,c_{s+1},\ldots,c_{2s}$ are its right
pole and left zero semi-residues.
\end{proof}
Using Lemma~\ref{prinrat}, we can give
the complete description of
 the class of Fuchsian systems, whose fundamental
solutions are rational matrix functions in general position,
in terms of the coefficients. The appropriate result,
formulated below, was  obtained by V.~E.~Katsnelson
in \cite{katze1}, \cite{katze}.
\begin{thrm}[Katsnelson] \label{t1.2}
Let $t_1,\ldots,t_{2s}$ be $2s$ distinct points in ${\mathbb C}$
and let $Q_1,\ldots,Q_{2s}\in gl(m,{\mathbb C}).$ Let $Y(x)$ be a
fundamental solution
  of the Fuchsian system
\begin{equation}\label{Furat1}
\der{Y}{x}=\sum_{j=1}^{2s}\frac{Q_{j}}{x-t_{j}}Y.\end{equation}
 Then $Y(x)$ is a rational matrix
function    in general position with the pole and zero sets
\begin{eqnarray}
{\cal P}&=&\Set{t_1,\ldots,t_s}, \\
{\cal Z}&=&\Set{t_{s+1},\ldots,t_{2s}} \end{eqnarray}
  if, and only if,  $Q_1,\ldots,Q_{2s}$ satisfy
  the following relations: \begin{eqnarray}
\label{regrat}\sum_{j=1}^{2s}Q_j&=&0,\\
\label{begcon} \rank{Q_{j}}&=&1,\quad j=1,\ldots,2s, \\
\label{reson} Q_{j}^2&=&\left\{ \begin{array}{l@{\quad}l} -
Q_{j},& j=1,\ldots,s, \\ Q_{j},&j=s+1,\ldots,2s, \end{array}
\right. \\ \label{sufconrat} Q_{j}R_{j}Q_{j}&=&\left\{
\begin{array}{l@{\quad}l} -R_{j}Q_{j},& j=1,\ldots,s, \\
Q_{j}R_{j},&j=s+1,\ldots,2s, \end{array} \right. \end{eqnarray}
where $R_1,\ldots,R_{2s}\in gl(m,{\mathbb C})$ are given by
\begin{equation}\label{endcon}
R_{j}=\sum_{j^{\prime}\not=j}\frac{Q_{j^{\prime}}}{t_{j}-t_{j^{\prime}}},\quad
j=1,\ldots,2s.
 \end{equation}
\end{thrm}
\begin{proof}
First of all, condition (\ref{regrat}) means that system~(\ref{Furat1})
is regular at $x=\infty.$  By definition, this condition is necessary
in order for $Y$ to be a rational matrix function in general position.
 Lemma~\ref{prinrat} implies that
conditions (\ref{begcon}), (\ref{reson}), which  mean that
for $j=1,\ldots,2s$ $Q_j$ has the Jordan form
\begin{equation}
J_j=\left\{ \begin{array}{l@{\quad}l} -I_1,& j=1,\ldots,s, \\
 I_1,&j=s+1,\ldots,2s, \end{array}
\right. \end{equation} are necessary as well. Thus we assume that
conditions (\ref{regrat}) - (\ref{reson}) are satisfied. In view
of
 Lemma~\ref{propfaclmm} in
Chapter~\ref{Schlesinger} and Lemma~\ref{prinrat}, we have to show
that for $j=1,\ldots,2s$ the linear system
\begin{equation}\label{holsysrat}
\der{H_j}{x}=\frac{Q_jH_j(x)-H_j(x)J_j}{x-
t_j}+\sum_{j^{\prime}\not=j}\frac{Q_{j^{\prime}}}
{x-t_{j^{\prime}}}H_j(x) \end{equation} has solution $H_j(x),$
holomorphic and invertible in a neighborhood of $t_j,$ if, and
only if, conditions (\ref{sufconrat}), (\ref{endcon}) are
satisfied. This can be done in the manner similar to the proof of
Lemma~\ref{Fufaclmm}. Let us fix $j,$ such that $1\leq j\leq s$
(in the case $s+1\leq j\leq 2s$ the proof is absolutely
analogous). Then we may assume, without loss of generality, that
\begin{equation}\label{assrat} Q_j=J_j=-I_1 \end{equation}
 and consider
the system \begin{equation}\label{holsysrat1}
\der{H_j}{x}=\frac{[H_j(x),I_1]}{x-
t_j}+\sum_{j^{\prime}\not=j}\frac{Q_{j^{\prime}}}
{x-t_{j^{\prime}}}H_j(x) \end{equation} with initial condition
\begin{equation}\label{ihsr} H_j(t_j)=I. \end{equation} Furthermore, we substitute
into (\ref{holsysrat1} ) the Taylor expansions
\begin{eqnarray}
\label{hjrat}
H_j(x)&=&I+\sum_{k=1}^{\infty}H_{j,k}(x-t_j)^k, \\
\sum_{j^{\prime}\not=j}\frac{Q_{j^{\prime}}} {x-
t_{j^{\prime}}}&=&R_j+\sum_{k=1}^{\infty}R_{j,k}(x-t_j)^k,
\end{eqnarray} where $R_j$ is given by (\ref{endcon}), and obtain
the recurrence relation for $H_{j,k}:$
\begin{equation}\label{recurrat}
(k+1)H_{j,k+1}-[H_{j,k+1},I_1]=R_jH_{j,k}+\sum_{k^{\prime}=1}^kR_{j,k^{\prime}}H_{j,k-
k^{\prime}}. \end{equation} Let us consider the linear operator
\begin{equation} \Psi_k:gl(m,{\mathbb C})\mapsto gl(m,{\mathbb C}),
\end{equation} depending on parameter $k\in{\mathbb Z}$ and
defined by \begin{equation} \Psi_k(X)=(k+1)X-[X,I_1].
\end{equation} We observe that
 $\Psi_k$ is invertible for $k>0$ and that
$X\in\Im{\Psi_0}$ if, and only if, \begin{equation} XI_1=I_1XI_1.
\end{equation} Hence, recurrence relation (\ref{recurrat}) is
solvable if, and only if, $R_j$ satisfies \begin{equation}
R_jI_1=I_1R_jI_1,
\end{equation} which is precisely the form taken by condition
(\ref{sufconrat}) under assumption (\ref{assrat}. It remains to
observe that if recurrence relation (\ref{recurrat}) is solvable,
then the convergence of series (\ref{hjrat}) and  the
invertibility of $H_j(x)$ in a neighborhood of $t_j$ follow from
the
 estimate
\begin{equation} \|\Psi_k^{-1}\|=O(\frac{1}{k}),\quad k\to +\infty,
\end{equation} on the norm of operator $\Psi_k^{-1}$ and from
initial condition (\ref{ihsr}), respectively.
\end{proof}
\begin{rem}
\label{solFurat}
Note that if conditions (\ref{regrat}) -
(\ref{endcon}) for system (\ref{Furat1}) are satisfied then
it is possible to obtain the explicit form of fundamental
solution $Y(x),$ satisfying the initial condition
(\ref{inirat}), directly from  $Q_1,\ldots,Q_{2s}.$
Indeed, in view of (\ref{begcon}), for $j=1,\ldots,2s$
$Q_j$ admits the factorization
\begin{equation}
Q_j=\left\{ \begin{array}{l@{\quad}l} c_jq_j,& j=1,\ldots,s, \\
 q_jb_j,&j=s+1,\ldots,2s, \end{array}
\right. \end{equation} where
$c_1,\ldots,c_s,q_{s+1},\ldots,q_{2s}$ and
$q_1,\ldots,q_{s},b_{s+1}, \ldots,b_{2s}$ are, respectively,
$m\times 1$ and $1\times m$ matrices. Furthermore,
Lemma~\ref{prinrat} implies that matrices
$c_1,\ldots,c_s;b_{s+1},\ldots,b_{2s}$ are the left pole and right
zero semi-residues of $Y.$ Hence, using (\ref{snp}), we can
construct from these semi-residues  zero-pole core matrix
$S_{{\cal Z}{\cal P}}$ of $Y(x)$ and determine,  according to
Theorem~\ref{t1.1}, the rest of the semi-residues of $Y(x).$
\end{rem}
\section{The Schlesinger System: Rational Solutions}
 Let us consider an isomonodromic deformation
of a Fuchsian system, whose fundamental solutions are rational
matrix functions in general position. Let set ${\cal S}$ be
defined by \begin{equation} {\cal S}={\mathbb
C}^{2s}\sm\bigcup_{j^{\prime}\not=j^{\prime\prime}}\Set{t=(t_1,\ldots,t_{2s}):t_{j^{\prime}
} =t_{j^{\prime\prime}}}. \end{equation}
 Let us fix
$t^0=(t_1^0,\ldots,t_{2s}^0)\in{\cal S}$ and $Q_1^0,\ldots,
Q_{2s}^0\in gl(m,{\mathbb C}),$ such that
\begin{eqnarray} \label{regrat0} \sum_{j=1}^{2s}Q_j^0&=&0, \\
\label{begcon0} \rank{Q_{j}^0}&=&1,\quad j=1,\ldots,2s, \\
\label{reson0} \left (Q_{j}^0\right )^2&=&\left\{
\begin{array}{l@{\quad}l} -Q_{j}^0,& j=1,\ldots,s, \\
Q_{j}^0,&j=s+1,\ldots,2s, \end{array} \right. \\
\label{sufconrat0} Q_{j}^0R_{j}^0Q_{j}^0&=&\left\{
\begin{array}{l@{\quad}l} -R_{j}^0Q_{j}^0,& j=1,\ldots,s, \\
Q_{j}^0R_{j}^0,&j=s+1,\ldots,2s, \end{array} \right.
\end{eqnarray} with $R_j^0$  given by \begin{equation}\label{endcon0}
R_{j}^0=\sum_{j^{\prime}\not=j}\frac{Q_{j^{\prime}}^0}{t_{j}^0-t_{j^{\prime}}^0}.
\end{equation} Let us consider solution $Y(x,t^0)$  of the Fuchsian system
\begin{equation}\label{Furat0} \der{Y}{x}=\sum_{j=1}^{2s}\frac{Q_{j}^0}{x-t_{j}^0}Y,
\end{equation} with initial condition \begin{equation} Y(\infty,t^0)=I. \end{equation} According to
Theorem~\ref{t1.2}, $Y(x,t^0)$ is a rational matrix function in
general position with the pole and zero sets
\begin{eqnarray}
\label{poleset0} {\cal P}^0&=&\Set{t_1^0,\ldots,t_s^0}, \\
\label{zeroset0} {\cal Z}^0&=&\Set{t_{s+1}^0,\ldots,t_{2s}^0}.
\end{eqnarray}
 In view of Remark~\ref{solFurat},
let us fix the factorization
\begin{equation}
Q_j^0=\left\{ \begin{array}{l@{\quad}l} c_j^0q_j^0,& j=1,\ldots,s, \\
 q_j^0b_j^0,&j=s+1,\ldots,2s, \end{array}
\right.
\end{equation}
where $c_1^0,\ldots,c_s^0,q_{s+1}^0,\ldots,q_{2s}^0$ and
$q_1^0,\ldots,q_{s}^0,
b_{s+1}^0,\ldots,b_{2s}^0$
are, respectively, $m\times 1$ and $1\times m$ matrices,
and construct  the left pole and right zero
semi-residual matrices of $Y(x,t^0):$
\begin{eqnarray}
C_{{\cal P}}^0&=&\left (c_1^0\cdots c_s^0\right ),
\\
B_{{\cal Z}}^0&=&\left (\begin{array}{c}b_{s+1}^0\\ \vdots \\
b_{2s}^0\end{array}\right ).
\end{eqnarray}
Furthermore, let ${\cal U}$ be a neighborhood of $t^0$ in ${\cal
S}$ and let $C_{{\cal P}}(t)$ and $B_{{\cal Z}}(t)$ be,
respectively, an $m\times s$ and an $s\times m$ matrix functions,
holomorphic in ${\cal U}$ and satisfying \begin{eqnarray}
\label{inlp} C_{{\cal P}}(t^0)&=&C_{{\cal P}}^0,
\\ \label{inrn} B_{{\cal Z}}(t^0)&=&B_{{\cal Z}}^0. \end{eqnarray} Let matrix
function $S_{{\cal Z}{\cal P}}(t)\in gl(s,{\cal O}({\cal U}))$
satisfy for each $t$ the equation \begin{equation}\label{lynpt}
A_{{\cal Z}}(t)S_{{\cal Z}{\cal P}}(t)-S_{{\cal Z}{\cal
P}}(t)A_{{\cal P}}(t)=
 B_{{\cal Z}}(t)C_{{\cal P}}(t),
\end{equation}
 where
\begin{eqnarray}
A_{{\cal P}}(t)&=&\diag{t_1,\ldots,t_s}, \label{apt} \\
A_{{\cal Z}}(t)&=&\diag{t_{s+1},\ldots,t_{2s}}. \label{ant}
\end{eqnarray} Then $S_{{\cal Z}{\cal P}}(t^0)$ is the zero-pole core matrix
of $Y(x,t^0).$ Hence, $\det{S_{{\cal Z}{\cal P}}(t^0)}\not=0$ and
we can assume that $S_{{\cal Z}{\cal P}}(t)\in GL(s,{\cal U}).$
According to Theorem~\ref{t1.1}, we can consider for each fixed
$t$ in ${\cal U}$ rational matrix function in general position
$Y(x,t),$ normalized by \begin{equation}\label{Furatdefini}
Y(\infty,t)\equiv I,
\end{equation} with  the pole and zero sets
\begin{eqnarray}
\label{polesett} {\cal P}&=&\Set{t_1,\ldots,t_s}, \\
\label{zerosett} {\cal
Z}&=&\Set{t_{s+1},\ldots,t_{2s}}\end{eqnarray} and left pole and
right zero semi-residual matrices $C_{{\cal P}}(t),B_{{\cal
Z}}(t).$ By Lemma~\ref{Fuchsrat}, $Y(x,t)$ satisfies the system
\begin{equation}\label{Furatdeform}
\dder{Y}{x}=\sum_{j=1}^{2s}\frac{Q_j(t)}{x-t_{j}}Y, \end{equation}
where
 $Q_1(t),\ldots,Q_{2s}(t)\in gl(m,{\cal O}({\cal U}))$ are given by
\begin{equation}\label{ratajt} Q_{j}(t)=\left\{ \begin{array}{l@{\quad}l} C_{{\cal
P}}(t)I_jS_{{\cal Z}{\cal P}}^{-1}(t)\left (t_{j}I-A_{{\cal
Z}}(t)\right )^{-1} B_{{\cal Z}}(t),&j=1,\ldots,s,\\
C_{{\cal P}}(t)\left (t_{j}I-A_{{\cal P}}(t)\right )^{- 1}S_{{\cal
Z}{\cal P}}^{-1}(t)I_{j-s} B_{{\cal Z}}(t),&j=s+1,\ldots,2s.
\end{array} \right. \end{equation} In particular, \begin{equation}\label{inischrat}
Q_j(t^0)=Q_j^0, \quad j=1,\dots,2s. \end{equation} Since the
monodromy representation of a rational matrix function in general
position is trivial, we conclude that system~(\ref{Furatdeform}),
thus constructed, gives in ${\cal U}$ an isomonodromic deformation
of Fuchsian system~(\ref{Ful0}). However, this deformation need
not be governed by the Schlesinger system. Indeed,  solution of
the latter must be completely determined at each $t$ by  initial
condition~(\ref{inischrat}), whereas in the construction above our
choice of $C_{{\cal P}}(t),B_{{\cal Z}}(t)$ is restricted  by
(\ref{inlp}), (\ref{inrn}) only at $t=t^0.$ Thus the conclusion of
Theorem~\ref{Schthrm} fails here. This does not, however,
contradict Theorem~\ref{Schthrm}, because (\ref{reson0}) means
that the non-resonance condition is violated.

The question remains,
how to find the solution of the Schlesinger system with initial
condition~(\ref{inischrat}), satisfying (\ref{regrat0}) -
(\ref{endcon0}). According to Theorem~\ref{schlthrm},
it is necessary and sufficient to find the isoprincipal deformation
of Fuchsian system (\ref{Furat0}), constructed from  initial
condition~(\ref{inischrat}). In view of Lemma~\ref{prinrat},
we can achieve this by modifying
the construction above so that the right pole and left
zero semi-residual matrices of $Y(x,t)$ are independent
of $t.$
 Thus, employing Theorem~\ref{t1.1},  we can prove the third
main result of this Thesis, obtained in joint work
with V.~E.~Katsnelson:
 \begin{thrm} \label{explicit} Let
$t^0=(t_1^0,\ldots,t_{2s}^0)\in{\cal S}$ and $Q_1^0,\ldots,
Q_{2s}^0\in gl(m,{\mathbb C})$ satisfy the relations
 \begin{eqnarray} \sum_{j=1}^{2s}Q_j^0&=&0, \\
\rank{Q_{j}^0}&=&1,\quad j=1,\ldots,2s, \\
\left (Q_{j}^0\right )^2&=&\left\{ \begin{array}{l@{\quad}l}
-Q_{j}^0,&
j=1,\ldots,s, \\ Q_{j}^0,&j=s+1,\ldots,2s, \end{array} \right. \\
Q_{j}^0R_{j}^0Q_{j}^0&=&\left\{
\begin{array}{l@{\quad}l} -R_{j}^0Q_{j}^0,& j=1,\ldots,s, \\
Q_{j}^0R_{j}^0,&j=s+1,\ldots,2s, \end{array} \right.
\end{eqnarray} with $R_j^0$  given by \begin{equation}
R_{j}^0=\sum_{j^{\prime}\not=j}\frac{Q_{j^{\prime}}^0}{t_{j}^0-t_{j^{\prime}}^0}.
\end{equation}
 Then solution $Q_1(t),\ldots,Q_{2s}(t)$
of the   Schlesinger system \begin{equation} \d Q_{j}
=\sum_{j^{\prime}\not=j} \left  [Q_{j^{\prime}}, Q_{j}\right
]\frac{\d t_{j}-\d t_j^{\prime}} {t_j-t_j^{\prime}},\quad
j=1,\ldots, 2s \end{equation} with the initial condition
\begin{equation} Q_j(t^0)=Q_j^0, \quad j=1,\dots,2s \end{equation}
can be constructed in the following way.
\begin{enumerate}
\item
  For $j=1,\ldots,2s$ factorize matrix $Q_j^0$ (of rank 1)
  in the form
  \begin{equation}
Q_j^0=\left\{ \begin{array}{l@{\quad}l} c_j^0q_j^0,& j=1,\ldots,s, \\
 q_j^0b_j^0,&j=s+1,\ldots,2s, \end{array}
\right.
\end{equation}
where $c_1^0,\ldots,c_s^0,q_{s+1}^0,\ldots,q_{2s}^0$ and
$q_1^0,\ldots,q_{s}^0,
b_{s+1}^0,\ldots,b_{2s}^0$
are, respectively, $m\times 1$ and $1\times m$ matrices.
   Form  $m\times s$
   matrix
  \begin{equation}
  C_{{\cal P}}^0=\left (c_1^0\cdots c_s^0\right )
  \end{equation}
  and $s\times m$ matrix
  \begin{equation}
  B_{{\cal Z}}^0 =\left (\begin{array}{c}b_{s+1}^0\\ \vdots \\ b_{2s}^0\end{array}\right ).
  \end{equation}
\item
  Define $s\times s$ matrix $S_{{\cal Z}{\cal P}}^0$ by the  equation
  \begin{equation}
  A_{{\cal Z}}^0S_{{\cal Z}{\cal P}}^0-S_{{\cal Z}{\cal P}}^0A_{{\cal P}}^0=
 B_{{\cal Z}}^0C_{{\cal P}}^0,
  \end{equation}
  where
  \begin{eqnarray}
  A_{{\cal P}}^0&=&\diag{t_1^0,\ldots,t_s^0}, \\
  A_{{\cal Z}}^0&=&\diag{t_{s+1}^0,\ldots,t_{2s}^0}. \end{eqnarray}
  That is, set
  \begin{equation}
  S_{{\cal Z}{\cal P}}^0=\left (\frac{b_{s+\alpha}^0c_{\beta}^0}
{t_{s+\alpha}^0-t_{\beta}^0}\right )_{\alpha,\beta=1}^s.
  \end{equation}
\item
  Determine the  matrix $S_{{\cal P}{\cal Z}}^0=\left (S_{{\cal Z}{\cal P}}^0\right )^{-1}$ and
  define $s\times m$
  matrix $B_{{\cal P}}^0$ and $m\times s$ matrix $C_{{\cal Z}}^0$ by
  \begin{eqnarray}
  B_{{\cal P}}^0&=&-S_{{\cal P}{\cal Z}}^0B_{{\cal Z}}^0,
  \\
  C_{{\cal Z}}^0 &=&
  C_{{\cal P}}^0S_{{\cal P}{\cal Z}}^0. \end{eqnarray}
\item
  Define
  rational $s\times s$ matrix function $ S_{{\cal P}{\cal Z}}(t)$ by the
   equation
  \begin{equation}
  A_{{\cal P}}(t)S_{{\cal P}{\cal Z}}(t)-S_{{\cal P}{\cal Z}}(t)A_{{\cal Z}}(t)=
  B_{{\cal P}}^0C_{{\cal Z}}^0,
  \end{equation}
  where
  \begin{eqnarray}
  A_{{\cal P}}(t)&=&\diag{t_1,\ldots,t_s}, \\
  A_{{\cal Z}}(t)&=&\diag{t_{s+1},\ldots,t_{2s}}. \end{eqnarray}
  That is, set
  \begin{equation}
   S_{{\cal P}{\cal Z}}(t)=\left (\frac{b_{\alpha}^0c_{s+\beta}^0}
{t_{\alpha}-t_{s+\beta}}\right )_{\alpha,\beta=1}^s,
  \end{equation}
  where $b_{\alpha}^0$ and $c_{s+\beta}^0$ are the $\alpha$-th row
of $B_{{\cal P}}^0$
  and the $\beta$-th column of $C_{{\cal Z}}^0,$ respectively.
\item
  Determine rational $s\times s$ matrix function
  $ S_{{\cal Z}{\cal P}}(t)= S_{{\cal P}{\cal Z}}^{-1}(t)$
  and define  rational $m\times s$ and $s\times m$
matrix functions $C_{{\cal P}}(t)$ and $B_{{\cal Z}}(t),$
respectively,
  by
  \begin{eqnarray}
  B_{{\cal Z}}(t)&=&-S_{{\cal Z}{\cal P}}(t)B_{{\cal P}}^0, \\ C_{{\cal P}}(t) &=&
  C_{{\cal Z}}^0S_{{\cal Z}{\cal P}}(t). \end{eqnarray}
\item
 Set
\begin{equation} Q_{j}(t)=\left\{ \begin{array}{l@{\ }l} C_{{\cal
P}}(t)I_jS_{{\cal P}{\cal Z}}(t)\left (t_{j}I-A_{{\cal
Z}}(t)\right )^{-1} B_{{\cal Z}}(t),&j=1,\ldots,s,\\
C_{{\cal P}}(t)\left (t_{j}I-A_{{\cal P}}(t)\right )^{- 1}S_{{\cal
P}{\cal Z}}(t)I_{j-s} B_{{\cal Z}}(t),&j=s+1,\ldots,2s,
\end{array} \right. \end{equation} where for $k=1,\ldots,s$ $I_k$ is the
diagonal $s\times s$ matrix \begin{equation}\begin{array}{rl}
I_k&=\left (\delta_{\alpha}^k\delta_{\beta}^k\right )_{\alpha,\beta=1}^s\\
&=\diag{0,\ldots,0,1,0,\ldots,0},
\end{array}
\end{equation}
 whose only non-zero element is in the $k$-th position on the diagonal.
 \end{enumerate}
\end{thrm}
\begin{rem}
\begin{enumerate}
\item In order to allow steps 3. and 5. of the construction in
Theorem~\ref{explicit}, matrix $S_{{\cal Z}{\cal P}}^0$ must be
invertible and the determinant of rational matrix function $
S_{{\cal P}{\cal Z}}(t)$ must not vanish identically.  We can see,
in view of Theorem~\ref{t1.1}, that this is indeed the case,
because, by our construction and by Theorem~\ref{t1.2}, $S_{{\cal
Z}{\cal P}}^0$ is the zero-pole core matrix of a rational matrix
function in general position and \begin{equation}
 S_{{\cal P}{\cal Z}}(t^0)=\left (S_{{\cal Z}{\cal P}}^0\right )^{-1}.
\end{equation}
\item There is a certain ambiguity in  the factorization of $Q_j^0$
at step 1. of
the construction in Theorem~\ref{explicit}.
It corresponds to the ambiguity in the definition of the semi-residual
matrices, described in Remark~\ref{amsemres}. However, this ambiguity
disappears in the definition of $Q_j(t)$ at step 6.
\item Note that solution $Q_1(t),\ldots,Q_{2s}(t)$ of the Schlesinger
system, constructed in Theorem~\ref{explicit}, is rational.
\end{enumerate}
\end{rem}
\bibliography{bibcor}
 \bibliographystyle{alpha}
  \end{document}